\documentclass[11pt]{amsart}

\title[Whittaker functions]
{Whittaker functions on orthogonal groups of odd degree}
\author{Taku Ishii}

\address{Chiba Institute of Technology, 2-1-1 Shibazono, Narashino, Chiba
275-0023, Japan}
\email{ishii.taku@it-chiba.ac.jp} 

\usepackage{amsmath}            \usepackage{amssymb}
\usepackage{amsopn}		\usepackage{amstext}
\usepackage{amsthm}		\usepackage{amscd}
\usepackage{amsgen}		\usepackage{latexsym}
\usepackage{enumerate}	
\usepackage{graphicx}

\newcommand{\newcom}{\newcommand}
\newcom{\renewcom}{\renewcommand}

\def\R{{\mathbb{R}}}     
\def\C{{\mathbb{C}}}     
\def\Z{{\mathbb{Z}}}     
\def\N{{\mathbb{N}}}     

\def\g{{\mathfrak{g}}}
\def\a{{\mathfrak{a}}}
\def\k{{\mathfrak{k}}}
\def\n{{\mathfrak{n}}}
\def\p{{\mathfrak{p}}}

\def\bm{{{\bf m}}}
\def\bk{{{\bf k}}}
\def\bl{{{\bf l}}}
\def\be{{{\bf e}}}

\theoremstyle{plain}
    
\newtheorem{thm}{Theorem}[section]
\newtheorem{lem}[thm]{Lemma}
\newtheorem{cor}[thm]{Corollary}
\newtheorem{prop}[thm]{Proposition}

\theoremstyle{definition}

\newtheorem{defn}[thm]{Definition}

\newcom{\bpf}{{\it Proof. }}  	
\newcom{\epf}{\hfill $\Box$}

\numberwithin{equation}{section}

\allowdisplaybreaks

\begin{document}
\maketitle
\begin{abstract}
We give explicit formulas for Whittaker functions for the class one principal
series representations of the orthogonal groups $ SO_{2n+1}(\R) $ of odd degree.
Our formulas are similar to the recursive formulas for Whittaker functions on
$SL_n(\R)$ given by Stade and the author \cite{ISt}.
Some parts of our results are announced in \cite{I3}.
\end{abstract}
\renewcommand{\thefootnote}{\fnsymbol{footnote}}
\section*{Introduction}

Special functions on a semisimple Lie group $G$ have been studied by many authors,
however, most of examples given in the literature are limited to the cases of 
rank one, 
and they are reduced to classical special functions
such as hypergeometric functions, Whittaker functions, Bessel functions and so on. 
As an effective example for higher rank case, 
we here give explicit formulas of certain special functions 
on orthogonal group $ SO_{2n+1}(\R) = SO_{n+1,n}(\R)$ of odd degree.

Let us briefly recall the spherical functions of Harish-Chandra \cite{HC}.
We fix a maximal compact subgroup $ K $ of $G$ and denote by
$ {\bf D}(G/K) $ the algebra of invariant differential operators on $ G/K $.
Let $ \chi_{\nu} : {\bf D}(G/K) \to \C $ be a homomorphism.
A smooth function $ f $ on $ G$ is called a spherical function if 
$ f(e)=1 $ and 
\begin{enumerate}[(1)]
\item $ f $ is bi-$K$ invariant, 
\item $ D f = \chi_{\nu}(D) f, \ \text{for all }  D \in {\bf D}(G/K). $ 
\end{enumerate}
Harish-Chandra gave the following integral representation for the spherical function
$$ 
 \psi_{\nu}(g) = \int_K a(kg)^{\nu+\rho} dk
$$
(see \S 1 for the notation), and obtained the expansion formula 
$$
 \psi_{\nu} = \sum_{w \in {\mathcal W}} c(w \nu) \Psi_{w \nu}, 
$$
where $ {\mathcal W} $ is the (small) Weyl group and 
$ \{ \Psi_{w \nu} \mid w \in {\mathcal W} \} $ is 
a basis of the solution of the system (2) and $ c(\nu) $ denotes
Harish-Chandra's $c$-function. 
The $c$-functions play many roles in harmonic analysis on $G$, for example 
they determine Plancherel measure for the spherical transform on $G$.

Our target in this paper is 
Whittaker function which is also a smooth function on $G$ satisfying the system (2).
Instead of the condition (1), we impose on  
$$
 f(ngk) = \eta(n)f(g), \ \text{for all } (n,g,k) \in N \times G \times K,
$$
where
$N$ is a maximal unipotent subgroup of $G$ and $ \eta $ a (nondegenerate)
unitary character of $N$. 

As is well known Whittaker functions appear in Fourier expansions of automorphic forms
(\cite{Sh}), and 
therefore they are indispensable in the various constructions of 
automorphic $L$-functions.
Jacquet \cite{Ja} introduced an integral representation of Whittaker function
of the form
\[
 J_{\nu,\eta}(g) = \int_{N} \eta^{-1}(n) \, a(w_0^{-1}ng)^{\nu+\rho} dn
\]
where $ w_0 $ is the longest element in $ {\mathcal W} $.
We refer to this Jacquet's Whittaker function as {\it class one Whittaker function}.
Inspired by the work of Harish-Chandra, Hashizume \cite{Ha} proved an 
expansion formula for the class one Whittaker function:
$$
 J_{\nu,\eta} = \sum_{w \in {\mathcal W}} c'(w  \nu) M_{w  \nu, \eta}.
$$
Here $ c'(\nu) $ is a product of $ c$-functions and certain ratio of gamma functions
and $ M_{\nu,\eta} $ is a power series solution (we call it {\it fundamental Whittaker function})
of the system (2) around the regular singularity.
We will recall the results of Hashizume in section 1.

Despite the development of the study of the expansion formulas, explicit formulas
of the spherical functions or Whittaker functions themselves seem to be still missing 
in most cases.
It is necessary to understand deeper the both sides of the expansion formula above
for serious applications to automorphic forms.

As for the fundamental Whittaker functions, 
recurrence relations characterizing the coefficients of them are easily given since they are 
essentially controlled by the Casimir operators. 
At the present these coefficients for 
$ SL_3(\R)$, $ SO_{5}(\R) $ and $ SL_4(\R) $ 
are known to be expressed in terms of 
(terminating) generalized hypergeometric
series $ _2F_1(1) $(=ratio of gamma functions by Gauss' formula), $ _3F_2(1)$
and $ _4F_3(1) $, respectively (see \cite{Bu}, \cite{I2}, \cite{St93}). 
But such a direction, that is, unit arguments of generalized  hypergeometric series 
do not seem be appropriate.
Recently Stade and the author \cite{ISt} reached a very satisfactory expression, 
which is a recursive relation between $ SL_{n-1}(\R) $ and $ SL_{n}(\R) $.
In section 2, a similar formula for $ SO_{2n+1}(\R)$ will be given.

Jacquet integrals are actually integral representations of 
class one Whittaker functions.
But it does not seem to be suitable form for applications to 
automorphic forms, such as computations of archimedean $L$-factors, or 
giving sharp estimates for Whittaker functions. 
Then we need to modify Jacquet integrals to more desirable expressions 
such as Mellin-Barnes type. 
In the case of $ SL_2(\R) $ the class one Whittaker function is essentially 
the modified 
$K$-Bessel function and it has the integral representations 
$$
 K_{\nu}(z) 
  = 2^{-1} \int_{0}^{\infty} 
     \exp \Bigl\{ -\frac{z}{2} \Bigl( t +\frac{1}{t} \Bigr) \Bigr\} t^{\nu} \frac{dt}{t}
  = \frac{2^{-2}}{2\pi\sqrt{-1}} \int_{-\sqrt{-1} \infty}^{\sqrt{-1} \infty} 
              \Gamma \Bigl( \frac{s+\nu}{2} \Bigr)
              \Gamma \Bigl( \frac{s-\nu}{2} \Bigr)
               \Bigl( \frac{z}{2} \Bigr)^{-s} ds.
$$

Extension to $SL_n(\R)$ was done by Stade \cite{St90}, \cite{St01}.
Starting from the Jacquet integral, he reached recursive formulas between 
class one Whittaker functions on 
$SL_n(\R)$ and $SL_{n-2}(\R)$ of both types above.
The formulas of Stade are very useful for 
an application, actually he computed certain archimedean zeta integrals 
in \cite{St01} and \cite{St02}.
Based on his formulas, Stade and the author \cite{ISt} find another 
recursive formula between $SL_n(\R) $ and $SL_{n-1}(\R)$ 
corresponding to that for fundamental Whittaker functions.

Our approach in this paper is different from \cite{ISt}, 
because we do not use the Jacquet integral.
Recall that the result of \cite{ISt} highly relies on the 
evaluation of the Jacquet integral 
which needs very complicated computation as shown in \cite{St90}.
Therefore we firstly try to guess a right formula by 
taking notice of an analogy between recursive formula for the 
fundamental Whittaker functions and the Mellin-Barnes integral representations of the 
class one Whittaker functions,
roughly speaking, the residue of the integrand of the Mellin-Barnes integral representation
gives the coefficient of the fundamental Whittaker functions.
Then we can arrive a conjectural form in view of the above integral representation
of $K$-Bessel functions (Theorem \ref{mainthm}).
This argument is of course heuristic and does not give a proof. 
In our proof in section 3, 
the expansion formulas of Hashizume and the recursive formula for the fundamental 
Whittaker functions obtained in section 2 play central roles. 
We remark that our idea of proof discussed here can be also 
applicable to the case of $SL_n(\R)$.

In section 4 we will observe
our recursive relation for $ SO_{2n+1}(\R) $ and $ SO_{2n-1}(\R) $
is similar to that for $ SL_{2n}(\R) $ and $ SL_{2n-2}(\R) $ in \cite{St90}. 
Following the argument of \cite{St01} we will compute the Mellin transforms of the class one 
Whittaker functions.

The author would like to thank Professor Takayuki Oda for his comment on the 
draft of the paper.

\section{preliminaries}
In this section we recall basic facts about
Whittaker functions for the class one principal series representations 
of $ SO_{2n+1}(\R) $. 
Our main reference is \cite{Ha}, which discusses for general
semisimple Lie groups.

\subsection{Group structures}

Let $ G = SO_{2n+1}(\R) $ be the special orthogonal group of 
degree $n$ with respect to an anti-diagonal matrix
$ \begin{pmatrix} 
  & & 1 \\ & \rotatebox[origin=c]{45}{$\cdots$} &  \\ 1 & & 
\end{pmatrix} $ of size $2n+1$.
Let $ \g = {\rm Lie}(G) = \k \oplus \p $ be a Cartan decomposition 
where $ \k $ and $ \p $ are $ +1 $ and $ -1 $ eigenspaces, respectively 
with respect to a Cartan involution $ \theta(X) = -{}^t X $ $ (X \in \g) $.
We take a maximal compact subgroup $K$ of $G$ by $ K = \exp \k \cong SO(n) \times SO(n+1) $.

Take a maximal abelian subalgebra 
$ \a = \{ {\rm diag} (t_1,\dotsc,t_n,0,-t_n,\dotsc,-t_1)  \mid t_i \in \R \} $ of $ \p $.
The restricted root system 
$ \Delta = \Delta(\g,\a) 
 = \{ \pm e_i \pm e_j, \, \pm e_k \mid  1 \le i<j \le n, \, 1 \le k \le n \}$ 
is of type $ B_n $, where $ e_i $ is a linear form on $ \a $ such that
$ e_i ( {\rm diag}(t_1,\dotsc,t_n,0,-t_n,\dotsc,-t_1) ) = t_i $. 
We take a positive system $ \Delta_+ $ and the simple system $ \Pi $ by 
$ \Delta_{+} = \{ e_i \pm e_j \mid  1 \le i<j \le n \} \cup
  \{ e_k \mid 1 \le k \le n \}$
and 
$ \Pi = \{ \alpha_i = e_i - e_{i+1} \mid 1 \le i \le n-1 \} \cup \{ \alpha_n = e_n \} $.
Denote by $ \g_{\alpha} $ the root space for $ \alpha \in \Delta $ and 
put $ \n = \sum_{ \alpha \in \Delta^+ } \g_{\alpha} $.
Then 
we have the Iwasawa decompositions $ \g = \n \oplus \a \oplus \k $ and 
$ G=NAK $ with 
$ N = \exp \n = \{ \mbox{upper triangular unipotent matrices in } G \} $
and 
$ A = \exp \a $.

\bigskip

Let $ \langle\ , \ \rangle $ be an inner product on 
the dual space $ \a^* $ of $ \a $
induced from the Killing form on $ \g$,
and we extend it to the complex dual $ \a_{\C}^{*} $.
Fix an element $ \nu $ of $ \a^*_{\C} $.
Since $ \a^*_{\C} \cong \C^n $, we can identify 
$ \nu $ with $  (\nu_1, \dotsc, \nu_n) \in \C^n $ via
$ \langle e_i, \nu \rangle /  \langle e_i, e_i  \rangle = 2\nu_i $. 
Let $ \rho_n = \frac{1}{2} \sum_{\alpha \in \Delta_+} \alpha $ be
the half sum of positive roots. Then
$$
 a^{\nu+\rho_n} = \exp (\nu+\rho_n)(\log a) 
  = \prod_{i=1}^{n} a_i^{2\nu_i+n-i+1/2} 
$$
for $ a = {\rm diag}(a_1,\dotsc, a_n, 1, a_n^{-1}, \dotsc, a_1^{-1}) \in A $.
We introduce a coordinate $ y = (y_1,\dotsc, y_n) $ on $ A $ by 
$$
  y_i = \frac{a_i}{a_{i+1}} \ (1 \le i \le n-1), \ \ y_n = a_n. 
$$  
Then 
$$
  a^{\nu+\rho_n} = y^{\nu+\rho_n}
   = \prod_{i=1}^n y_i^{2(\nu_1 + \dotsb + \nu_i) + i(n-i/2) }.
$$

We denote by 
$ \mathcal{W}_n = \langle w_i \mid 1 \le i \le n \rangle
 \cong  \mathfrak{S}_n \times (\Z/2\Z)^{n} $ 
the Weyl group of $ \Delta $.
Here $ w_i $ is the simple reflection with respect to the 
simple root $ \alpha_i $.
We note the action of $ w_i $ on $ \nu = (\nu_1,\dotsc,\nu_n) \in \a_{\C}^* $:
$ w_i \nu = (\nu_1, \dotsc, \nu_{i-1}, \nu_{i+1}, \nu_i, \nu_{i+2}, \dotsc, \nu_n) $ 
for $ 1 \le i \le n-1 $ and 
$ w_n \nu = (\nu_1,\dotsc,\nu_{n-1},-\nu_n) $.

A nondegenerate unitary character $ \eta $ of $N$ is of the form   
$$
 \eta \bigl( (u_{i,j}) \bigr) 
 = \exp \Bigl( 2\pi \sqrt{-1} \sum_{i=1}^{n} \eta_i u_{i,i+1}  \Bigr),
\ \ \ (u_{i,j}) \in N 
$$
with nonzero real numbers $ \eta_i $ $(1 \le i \le n )$.
We remark that in the notation of \cite{Ha}, 
$ |\eta_{\alpha_i}| = \sqrt{2} \pi \eta_i / \sqrt{2n-1} $.  

\bigskip 

Let $U(\g_{\C})$ and $U(\a_{\C})$ be the universal enveloping algebras 
of $\g_{\C} $ and $\a_{\C}$, 
the complexifications of $\g$ and $\a$, respectively.
Set
\[
  U(\g_{\C})^K = \{ X \in U(\g_{\C}) \mid \mathrm{Ad}(k) X = X \mbox{ for all } k \in K \}.
\]
Let $p$ be the projection $ U(\g_{\C}) \to U(\a_{\C}) $ along the decomposition
\[
  U(\g_{\C}) = U(\a_{\C}) \oplus ( \mathfrak{n} \>\! U(\g_{\C}) + U(\g_{\C}) \>\! \k).
\]
Define an automorphism $ \gamma $ of $ U(\a_{\C}) $ by 
$ \gamma(H) = H + \rho_n(H) $
for $H \in \a_{\C}$. 
For the linear form $ \nu $ above we define an algebra homomorphism 
$ \chi_{\nu} : U(\g_{\C})^K \to \C $ by
\[
   \chi_{\nu} (z) = \nu( \gamma \circ p (z) ), \ \ \ z \in U(\g_{\C})^K.
\]
Note that $ \chi_{\nu} $ is trivial on $ U(\g)^K \cap U(\g) \>\! \k$.

\begin{defn}
Under the notation above we denote by $ {\rm Wh}(\nu,\eta) $ 
the space of smooth functions $ f : G \to \C $ satisfying 
\begin{itemize}
\item $f(ngk) = \eta(n) f(g) $ for all $ n \in N $, $ g \in G $ and $ k \in K $,
\item $Zf=\chi_{\nu}(Z)f $ for all $ Z \in U(\g_{\C})^K $.
\end{itemize}
\end{defn}


The Iwasawa decomposition $ G=NAK$ implies that any function
$ f $ in the space $ {\rm Wh}(\nu,\eta)$ 
is determined by its restriction $ f|_A $ to $ A $, which we call the 
{\it  radial part} of $f$.

\bigskip 

We mention a relation with the principal series representation of $G$.
Let $ M $ be the centralizer of $K$ in $A$ and $ P_{\rm min} = MAN $ 
a minimal parabolic subgroup of $G$.
The induced representation
$ \pi_{\nu} = {\rm Ind}_{P_{\rm min}}^G(1_M \otimes \exp(\nu+\rho_n) \otimes 1_N) $
is called the {\it class one principal series representation} of $G$. 
Consider an intertwining space $ {\rm Hom}_{G}(\pi_{\nu}, {\rm Ind}_N^G(\eta)) $, where
$ {\rm Ind}_N^G(\eta) 
 = \{ f \in C^{\infty}(G,\C) \mid f(ng) = \eta(n) f(g)  \mbox{ for all } (n,g) \in N \times G\}. $ 
Then our target space $  {\rm Wh}(\nu,\eta) $ 
can be thought as a realization of $ \pi_{\nu} $ in the induced module $ {\rm Ind}_N^G(\eta)  $, 
that is, for a nonzero intertwiner 
$ \Phi \in {\rm Hom}_{G}(\pi_{\nu}, {\rm Ind}_N^G(\eta)) $, 
and the spherical vector $ v \in \pi_{\nu} $, $ \Phi(v) $ becomes a nonzero element
in $ {\rm Wh}(\nu, \eta) $.

\subsection{Fundamental Whittaker functions}

Let us recall Hashizume's construction of a basis of the space $ \mathrm{Wh}(\nu, \eta) $
in our situation (\cite[\S 4]{Ha}).

For a set of nonnegative integers $ \bm = (m_1, \dotsc, m_n) $, 
we determine complex numbers $ c_{n,\bm}(\nu) $ by 
the initial condition $ c_{n,(0,\dotsc,0)} (\nu) = 1 $ and the recurrence relation
\begin{equation} \label{rec}
 \displaystyle q_n ( \bm, \nu ) c_{n,\bm} (\nu)
 = \sum_{i=1}^{n-1} c_{n,\bm-\be_i}(\nu) + \frac12 c_{n,\bm-\be_n}(\nu),
\end{equation}
where $ \be_i $ $ (1\le i \le n) $ is the $i$-th standard basis in $\R^n$
and 
$ q_n $ is defined by
\begin{align*}
 q_n(\bm, \nu)
& \equiv q_n \bigl( (m_1, \dotsc, m_n), (\nu_1, \dotsc, \nu_n) \bigr)
\\
& := \sum_{i=1}^{n-1} m_i^2 + \frac{1}{2} m_n^2 
   - \sum_{i=1}^{n-1} m_i m_{i+1} 
   + \sum_{i=1}^{n-1} (\nu_i-\nu_{i+1}) m_i
   + \nu_n m_n. 
\end{align*}
Hereafter we sometimes use the same symbol $ \be_i $ 
for the $i$-th standard basis in $\R^{n-1}$.
If $ q_n(\bm,\nu) $ does not vanish for all nonzero $ \bm $, 
we can uniquely determine $ c_{n,\bm}(\nu) $.

\begin{defn}
Define a power series 
$ M_{\nu,\eta}^{n}(y) = y^{\rho_n} \widetilde{M}_{\nu,\eta}^n(y) $ on $ A $ by
\begin{align*}
 \widetilde{M}_{\nu,\eta}^n(y) 
&:=  \sum_{m_1, \dotsc, m_n= 0}^{\infty} c_{n,(m_1,\dotsc, m_n)}(\nu)
    \prod_{i=1}^{n-1} (\pi \eta_i y_i)^{2 (m_i + \nu_1 + \cdots + \nu_i)}
    \cdot (\sqrt{2} \pi \eta_n y_n)^{2 (m_n + \nu_1 + \cdots + \nu_n)}
\end{align*}
and extend it to a function on $G$ by
$$
 M_{\nu,\eta}^n(g) = \eta(n(g)) M_{\nu,\eta}^n(a(g)),
$$
where we denote by $$ g=n(g) a(g) k(g), \ \ n(g) \in N, \, a(g) \in A, \, k(g) \in K $$
the Iwasawa decomposition of $ g $.
We call the function $ M_{\nu,\eta}^n $ the {\it fundamental Whittaker function} on $G$. 
\end{defn}

It is known that (\cite[Lemma 4.6]{Ha}) the power series $ M_{\nu,\eta}^{n}(y) $ 
converges absolutely and uniformly as functions of $ y \in A $ and 
$ \nu \in  \a_{\C}^* $ (cf. Lemma \ref{conv}).

\begin{defn}
An element $ \nu $ of $ \a_{\C}^{*} $ is called {\it regular} 
if the following two conditions are satisfied.
\begin{itemize}
\item 
$ q_n(\bm,w \nu) \neq 0  $
for all $ \bm  \neq (0,\dotsc,0) $ and $ w \in {\mathcal W}_n$,
\item
$ w\nu - w'\nu \notin \{ \sum_{i=1}^{n} m_{i} \alpha_i \mid m_{i} \in \Z \} $
for all pairs $ (w,w') $ in $ {\mathcal W}_n $ with $ w \neq w' $.
\end{itemize} 
We denote by $\,\!^{\prime} \!\!\> \a_{\C}^* $ 
the subset of regular elements in $\a_{\C}^{*}$.
\end{defn}

Hashizume proved the following.

\begin{thm} $($\cite[Theorem 5.4]{Ha}$)$
If $ \nu $ is a regular element, then the set
$$
  \{ M_{w\nu,\eta}^n  \mid  w \in {\mathcal W}_n \}
$$
forms a basis of $ {\rm Wh}(\nu,\eta) $.

\end{thm}


\subsection{Jacquet integral} 

It is known that the subspace $ {\rm Wh}(\nu,\eta)^{\rm mod} $
of moderate growth functions (\cite{Wa})
of $ {\rm Wh}(\nu,\eta) $ is at most one dimensional 
(\cite{Sh}, \cite{Ko}). 
Jacquet \cite{Ja} introduced an integral representation of the unique
element in $ {\rm Wh}(\nu,\eta)^{\rm mod} $:
$$
 J_{\nu,\eta}(g) = \int_N  \eta^{-1}(n) \,a(w_0^{-1} n g)^{\nu + \rho_n} dn,
$$
where
$ w_0 $ is the longest element in ${\mathcal W}_n$
and
$ dn $ is a normalized Haar measure on $N$ as in \cite[\S 1]{Ha}.


\begin{prop} $($\cite[Proposition 4.2]{Ha0}$)$
Let $D$ be a subset of $\a^{*}_{\C}$ defined by
$$
  D = \{ \nu \in \a_{\C}^* \mid \mathrm{Re}(\nu_i \pm \nu_j ) > 0 \ (1\le i<j \le n),
          \ \mathrm{Re}(\nu_i)>0  \ (1 \le i \le n) \}.
$$
Then the Jacquet integral $ J_{\nu,\eta} $
converges absolutely and uniformly on $ (g,\nu) \in G \times D $ and gives 
a holomorphic function on $\nu \in D$. 
Moreover, as a function on $ \nu $, 
$ J_{\nu,\eta} $ can be continued to an entire function on $ \a_{\C}^*$
and satisfies a functional equation
$$
 J_{\nu,\eta}(g) = \gamma(w, \nu, \eta)  J_{w\nu,\eta}(g) 
$$
for $ w \in {\mathcal W}_n $ and $ g \in G$.
Here $ \gamma(w, \nu,\eta) $ is defined as follows. 
For the simple reflection $w_i $, 
\begin{align*}
 \gamma(w_i, \nu,\eta) 
= 
\begin{cases}
  ( \pi \eta_i )^{2(\nu_i-\nu_{i+1})} 
  \dfrac{ \Gamma( -\nu_{i}-\nu_{i+1} + 1/2 ) }{ \Gamma( \nu_{i} + \nu_{i+1} + 1/2 ) }
 & \mbox{ if } 1 \le i \le n-1,  \\
 & \vspace{-3mm} \\
  ( \sqrt{2} \pi \eta_n )^{4\nu_n} 
  \dfrac{ \Gamma( -2 \nu_{n} + 1/2 ) }{ \Gamma( 2 \nu_{n} + 1/2 ) }
 & \mbox{ if } i = n,
\end{cases}
\end{align*}
and for $ w \in {\mathcal W}_n $ 
with $ l(w_{i} w) =l(w) +1 $, 
\[
  \gamma(w_{i} w , \nu,\eta) 
    = \gamma(w, \nu,\eta) \gamma(w_i,  w\nu, \eta),
\]         
where $ l(w) $ means the length of $w$.
\end{prop}

\begin{defn}
We call the Jacquet integral $ J_{\nu,\eta} $ (and its constant multiple) 
the {\it class one Whittaker function} on $G$.
\end{defn}

As in the way of Harish-Chandra, Hashizume expressed 
the class one Whittaker function $ J_{\nu,\eta} $ as a linear combination of 
the fundamental Whittaker functions $ M_{w \nu, \eta}^n $ $ (w \in {\mathcal W}_n) $.

\begin{thm} $($\cite[Theorem 7.8]{Ha}$)$
If $ \nu  $ is a regular element, then we have 
\begin{align*}
 J_{\nu,\eta}(g) 
& = \sum_{w \in {\mathcal W}_n} 
 \gamma(w_0 w, \nu, \eta) c(w_0 w \nu) M_{w\nu, \eta}^n(g),
\end{align*}
where
\begin{align*}
 c(\nu)
  &= \int_{N} a(w_0^{-1}n)^{\nu+\rho_n} dn
\\ 
 & = 2^{n/2} \{(4n-2) \pi\}^{n(n+1)/4}  
   \prod_{1 \le i<j \le n} 
    \frac{ \Gamma(\nu_i-\nu_j) \Gamma(\nu_i+\nu_j) }
         { \Gamma(\nu_i-\nu_j+1/2)  \Gamma(\nu_i+\nu_j+1/2) }
   \prod_{1 \le i \le n }
    \frac{ \Gamma(2\nu_i) }{ \Gamma(2\nu_i+1/2) }
\end{align*}
is the Harish-Chandra $c$-function on $G$.
Equivalently, if we put 
\begin{align*}
 W_{\nu,\eta}^n(g) 
& := 2^{-n/2} \{(4n-2)\pi\}^{-n(n+1)/4}
  \prod_{1 \le i \le n} \pi^{-2(\nu_1+ \dotsb + \nu_i)}
\\
& \quad \cdot 
  \prod_{1 \le i<j \le n} \Gamma(\nu_i-\nu_j+1/2) \Gamma(\nu_i+\nu_j+1/2) 
  \prod_{1 \le i \le n} \Gamma(2\nu_i +1/2)
  \cdot J_{\nu,\eta}(g), 
\end{align*}
then we have 
$$
 W_{\nu,\eta}^n(g) = \sum_{w \in {\mathcal W}_n} 
  w \Biggl[ \prod_{1 \le i<j \le n} \Gamma(-\nu_i+\nu_j) \Gamma(-\nu_i-\nu_j)
           \prod_{1 \le i \le n} \Gamma( -2\nu_i ) 
          \cdot M_{\nu, \eta}^n(g) \Biggr].
$$
\end{thm}

From now on we assume $ \eta_i = 1 $ for $ 1 \le i \le n-1 $ 
and $ \eta_n = 1/\sqrt{2} $ for simplicity and 
denote by $ M_{\nu}^n = M_{\nu,\eta}^n $,
$ W_{\nu}^n = W_{\nu,\eta}^n $ omitting the symbol $\eta$.

\section{Explicit formulas for fundamental Whittaker functions}
In this section we solve the recurrence relation (\ref{rec})
to find an explicit formula for fundamental Whittaker function.
Similar formulas for $SL_n(\R)$ are given in \cite{I1} and \cite{ISt}.
We use the Pochhammer symbol $ (a)_n = \Gamma(a+n)/ \Gamma(a) $ for $ n \in \Z $.

\begin{thm} \label{fund}
For $ \nu = (\nu_1,\dotsc,\nu_n)  \in \,\!^{\prime} \!\!\> \a_{\C}^* $,
put $ \tilde{\nu} = (\nu_1, \dotsc, \nu_{n-1}) $. 
Then we have
$ c_{1,m_1}(\nu) = 1/m_1!(2\nu_1+1)_{m_1} $ and 
\begin{align} \label{fundformula}
\begin{split}
 c_{n,(m_1, \dotsc, m_n)}(\nu)
& = 
 \sum_{ \scriptstyle \{l_1, \dotsc, l_{n-1} \} \atop
        \scriptstyle \{k_1, \dotsc, k_{n-1} \} }
  \frac{ c_{ n-1, (k_1, \dotsc, k_{n-1}) } (\tilde{\nu}) }
  { \prod_{i=1}^{n-1} (m_i-l_i)! \cdot (m_n-k_{n-1})! \, 
    \prod_{i=1}^{n-1} (l_i-k_i)! }
\\
& \qquad \qquad 
   \cdot \frac{1}{
    \prod_{i=1}^{n} ( \nu_i+\nu_n +1)_{m_i-l_{i-1}}
    \prod_{i=1}^{n-1}( \nu_i-\nu_n +1)_{l_i-k_{i-1}} },
\end{split}
\end{align} 
where 
the indexing sets $ \{ k_i \}$ and $ \{ l_i \} $ are the sets of nonnegative integers
satisfying
$$
   0 \le k_i \le l_i \le m_i \ (1 \le i \le n-1), 
 \ \  0 \le k_{n-1} \le m_n 
$$
and we promise $ k_0 = l_0 = 0 $.
\end{thm}

\bpf 
The idea of proof is similar to \cite{I1} and \cite{ISt} and 
we will check the right hand side of (\ref{fundformula}) satisfies 
(\ref{rec}) by dividing into two steps.

For a set of nonnegative integers $ \bl = (l_1,\dotsc,l_n) $ with
$ 0 \leq l_{n} \leq l_{n-1} $, define 
$ b_{\bl}(\nu)  \equiv b_{(l_1, \dotsc, l_n)}(\nu) $ by 
\begin{align*}
 b_{\bl}(\nu) 
& := (-1)^{l_n} 
  \sum_{ \{k_1,\dotsc, k_{n-2} \} }
  \frac{ c_{n-1, (k_1, \dotsc, k_{n-2}, l_n)} (\tilde{\nu}) }
  { \prod_{i=1}^{n-2} (l_i-k_i)!  \cdot  (l_{n-1}-l_n)! 
    \prod_{i=1}^{n-1} ( \nu_i-\nu_n + 1 )_{l_i-k_{i-1}} },
\end{align*}
where $ \{k_1,\dotsc, k_{n-2} \} $ means $ 0 \le k_i \le l_i $ for $ 1 \le i \le n-2 $.
The term $(-1)^{l_n} $ is included for our later convenience. 
Then the formula (\ref{fundformula}) is equivalent to 
\begin{align} \label{two}
 c_{n,(m_1,\dotsc, m_n)} (\nu) 
 = \sum_{\{ l_1, \dotsc, l_n \}  }
  \frac{ (-1)^{l_n} b_{(l_1, \dotsc, l_n)} (\nu) }
  { \prod_{i=1}^{n} (m_i-l_i)! 
    \prod_{i=1}^{n} ( \nu_i+\nu_n + 1 )_{m_i-l_{i-1}} },
\end{align}
where $ \{ l_1, \dotsc, l_n \} $ means $ 0 \le l_i \le m_i $ for $  1 \le i \le n $.

\smallskip

We first show that $ b_{\bl}(\nu) $ satisfies the recurrence relation
\begin{align} \label{recb}
  q_n ( \bl, \nu )
  b_{\bl}(\nu) 
= \sum_{i=1}^{n-1} b_{\bl-\be_i} (\nu) 
  + \frac{1}{2}(-l_{n-1}+l_n-1) b_{\bl-\be_n} (\nu).
\end{align}
Set
$$ P_{\bl, \bk}(\nu) = 
  \prod_{i=1}^{n-2} (l_i-k_i)! \cdot (l_{n-1}-l_n)! 
  \prod_{i=1}^{n-1} ( \nu_i-\nu_n + 1)_{l_i-k_{i-1}},
$$
the denominator of the summand in $ b_{\bf l}(\nu) $.
The key identity is
\begin{equation} \label{key1}
\begin{split}
& \sum_{i=1}^{n-1} \frac{P_{\bl, \bk}(\nu)}{ P_{\bl-\be_i, \bk}(\nu) } 
 - \sum_{i=1}^{n-2} \frac{P_{\bl,\bk}(\nu)}{ P_{\bl, \bk+\be_i}(\nu) }
= q_n \bigl( (l_1, \dotsc, l_n),  \nu \bigr)
  - q_{n-1} \bigl( (k_1, \dotsc, k_{n-2}, l_n), \tilde{\nu} \bigr).
\end{split}
\end{equation}
This is an easy algebra since the left hand side of (\ref{key1})
can be written as
\begin{equation*}
\begin{split}
&  \sum_{i=1}^{n-2} (l_i-k_i)(l_i-k_{i-1} + \nu_i-\nu_n)
  - \sum_{i=1}^{n-2} (l_i-k_i)(l_{i+1}-k_i+\nu_{i+1}-\nu_n)
\\
& + (l_{n-1}-l_n)(l_{n-1}-k_{n-2}+ \nu_{n-1}-\nu_n).
\end{split}
\end{equation*}
Let us compute the right hand side of (\ref{recb}).
Since
\begin{align*}
 b_{\bl-\be_i}(\nu)
&= 
 (-1)^{l_n} 
 \sum_{ \{ k_1, \dotsc, k_{n-2} \}  } 
 \frac{ P_{\bl, \bk}(\nu) }{ P_{\bl-\be_{i}, \bk}(\nu) }
 \cdot
 \frac{ c_{n-1, (k_1,\dotsc,k_{n-2},l_n)} ( \tilde{\nu} ) }
      { P_{\bl, \bk}(\nu) }
\end{align*} 
for $  1 \le i \le n-1 $    
and 
\begin{align*}
& 
 \frac{1}{2} (-l_{n-1}+l_n-1 ) b_{\bl-\be_n}(\nu)
\\
&= \frac{1}{2} (-l_{n-1}+l_n-1 ) 
   \cdot (-1)^{{l_n}-1}
   \sum_{ \{ k_1, \dotsc, k_{n-2} \}  } 
   \frac{ P_{\bl, \bk}(\nu) }{ P_{\bl-\be_{n}, \bk}(\nu) }
   \cdot
   \frac{ c_{n-1, (k_1,\dotsc,k_{n-2},l_n-1)} ( \tilde{\nu} ) }
       { P_{\bl, \bk}(\nu) }
\\
&= \frac{1}{2} (-l_{n-1}+l_n-1 ) 
   \cdot (-1)^{{l_n}-1}
   \sum_{ \{ k_1, \dotsc, k_{n-2} \}  } 
   \frac{1}{l_{n-1}-l_n+1}
   \cdot
   \frac{ c_{n-1, (k_1,\dotsc,k_{n-2},l_n-1)} ( \tilde{\nu} ) }
       { P_{\bl, \bk}(\nu) }
\\
&= (-1)^{l_n}
   \sum_{ \{ k_1, \dotsc, k_{n-2} \}  } 
   \frac{ \frac12  c_{n-1, (k_1,\dotsc,k_{n-2},l_n-1)} ( \tilde{\nu} ) }
        { P_{\bl, \bk}(\nu) },
\end{align*}
the identity (\ref{key1}) implies 
that the right hand side of (\ref{recb}) can be written 
as a sum of the following four terms:
\begin{align} \label{sum1}
 (-1)^{l_n}
 \sum_{i=1}^{n-2} 
 \sum_{ \{ k_1, \dotsc, k_{n-2} \} }
  \frac{P_{\bl,\bk}(\nu)}{ P_{\bl, \bk+\be_i}(\nu) }
  \cdot
  \frac{ c_{n-1 (k_1,\dotsc,k_{n-2},l_n)} ( \tilde{\nu} ) }
       { P_{\bl, \bk}(\nu) },
\end{align}
\begin{align} \label{sum2}
 (-1)^{l_n}
   \sum_{ \{ k_1, \dotsc, k_{n-2} \}  } 
   \frac{ \frac12  c_{n-1, (k_1,\dotsc,k_{n-2},l_n-1)} ( \tilde{\nu} ) }
        { P_{\bl, \bk}(\nu) },
\end{align}
\begin{align} \label{sum3}
 q_n \bigl( (l_1, \dotsc, l_n),  \nu \bigr)
 \cdot (-1)^{l_n}
   \sum_{ \{ k_1, \dotsc, k_{n-2} \}  } 
   \frac{ c_{n-1, (k_1,\dotsc,k_{n-2},l_n)} ( \tilde{\nu} ) }
        { P_{\bl, \bk}(\nu) }
  = q_n(\bl,\nu) b_{\bl}(\nu),
\end{align}
and 
\begin{align} \label{sum4}
 - q_{n-1} \bigl( (k_1, \dotsc, k_{n-2}, l_n), \tilde{\nu} )
  \cdot (-1)^{l_n} 
  \sum_{ \{ k_1, \dotsc, k_{n-2} \}  } 
   \frac{ c_{ n-1, (k_1,\dotsc,k_{n-2},l_n)} ( \tilde{\nu} ) }
        { P_{\bl, \bk}(\nu) }.
\end{align}
In (\ref{sum1}), we substitute $ k_i \to  k_i-1 $ to rewrite 
\begin{align*}
 (-1)^{l_n}
 \sum_{i=1}^{n-2} 
 \sum_{ \{ k_1, \dotsc, k_{n-2} \} }
  \frac{ c_{n-1, (k_1,\dotsc,k_{n-2},l_n) - \be_i } ( \tilde{\nu} ) }
       { P_{\bl, \bk}(\nu) }.
\end{align*}
Thus in view of the recurrence relation for 
$ c_{n-1, (k_1,\dotsc,k_{n-2},l_n) } ( \tilde{\nu} ) $, we find 
$ (\ref{sum1}) + (\ref{sum2}) + (\ref{sum4}) = 0 $ and finish the proof of 
(\ref{recb}).

\medskip 

In the next step we prove the right hand side of (\ref{two}) satisfies the recurrence 
relation (\ref{rec}) for $ c_{n,\bm}(\nu) $.
As in the first step, if we put 
$$
  Q_{\bm, \bl}(\nu) 
 = \prod_{i=1}^{n}(m_i-l_i)! 
    \prod_{i=1}^{n}(\nu_i + \nu_n +1)_{m_i-l_{i-1}}
$$
then the identity
\begin{align*} 
& \sum_{i=1}^{n-1}
   \frac{ Q_{\bm, \bl}(\nu) }{ Q_{\bm-\be_i, \bl}(\nu) }
+\frac12
   \frac{ Q_{\bm, \bl}(\nu) }{ Q_{\bm-\be_n, \bl}(\nu) }
-\sum_{i=1}^{n-1}
   \frac{ Q_{\bm, \bl}(\nu) }{ Q_{\bm, \bl+\be_i}(\nu) }
-\frac12 (l_{n-1}-l_n)(m_n-l_n)
\\
& =
  q_n(\bm, \nu) - q_n(\bl, \nu)
\end{align*}
holds. By means of (\ref{recb}) our claim follows and thus
we complete the proof of Theorem \ref{fund}.
\epf


\section{Explicit formulas for class one Whittaker functions}

In this section we will show a recursive integral representation of 
the class one Whittaker function.

\begin{thm} \label{mainthm} 
For $ \nu = (\nu_1, \dotsc,\nu_n) \in \a_{\C}^*$ and $ y = (y_1,\dotsc,y_n) \in A $,
we inductively define a function 
$  \widetilde{W}_{\nu}^n(y)$ on $A$ by
\begin{equation} \label{intrep1}
\begin{split} 
 \widetilde{W}_{\nu}^{n}(y) 
& := \int_{(\R^+)^{n}}
    \int_{(\R^+)^{n-1}}
    \prod_{i=1}^{n} \exp \biggl\{ -(\pi y_i)^2 t_i - \frac{1}{t_i} \biggr\}
    \prod_{i=1}^{n} (\pi y_i)^{2\nu_n}
\\
& \qquad \cdot  
  \prod_{i=1}^{n-1} 
     \exp \biggl\{ -(\pi y_i)^2 \frac{t_i}{t_{i+1}} u_i - \frac{1}{u_i} \biggr\}
  \cdot t_1^{2\nu_n} \prod_{i=1}^{n-1} (t_{i+1} u_i)^{\nu_n}
\\
& \qquad \cdot 
  \widetilde{W}_{\tilde{\nu} }^{n-1} 
  \Biggl( y_2 \sqrt{\frac{t_2 u_2}{t_3 u_1}}, \dotsc, 
          y_{n-1} \sqrt{\frac{t_{n-1} u_{n-1}}{ t_n u_{n-2}}},
          y_n \sqrt{ \frac{t_n}{u_{n-1}} } \Biggr) 
  \prod_{i=1}^{n-1} \frac{du_i}{u_i} \prod_{i=1}^{n} \frac{dt_i}{t_i},
\end{split}
\end{equation}
and $ \widetilde{W}_{\nu}^{1}(y) = \widetilde{W}_{\nu_1}^{1}(y_1) 
= 2 K_{2\nu_1}(2\pi y_1) $.
Here $ \tilde{\nu} = (\nu_1, \dotsc, \nu_{n-1}) $.
Then we have
\begin{align} \label{exp}
 \widetilde{W}_{\nu}^{n}(y) 
= \sum_{w \in {\mathcal W}_n}
  w \Bigl[ 
    \Gamma_n(\nu)
    \cdot  \widetilde{M}_{\nu}^n (y) \Bigr]                       
\end{align}  
with 
$$ \Gamma_n(\nu) := \prod_{1 \le i<j \le n} 
      \Gamma ( -\nu_i-\nu_j )
      \Gamma ( -\nu_i+\nu_j )
    \prod_{1 \le i \le n} \Gamma(-2 \nu_i)
$$
and thus $ W_{\nu}^n(y) = y^{\rho_n} \widetilde{W}_{\nu}^n(y) $.
\end{thm}
\medskip 

We illustrate the outline of the proof of the expansion formula (\ref{exp}).
It is done by induction on $n$ and 
as in the proof of Theorem \ref{fund}, it consists of two steps. 
For $ x = (x_1, \dotsc, x_n) \in (\R^+)^{n}$, let us define a function 
$ V_{\nu} (x)  = V_{(\nu_1,\dotsc,\nu_n)}(x_1,\dotsc,x_n) $ by 
\begin{equation} \label{defV}
\begin{split}
 V_{\nu} (x) 
& :=  \int_{(\R^+)^{n-1}}
  \prod_{i=1}^{n-1} \exp \biggl\{ -(\pi x_i)^2 u_i- \frac{1}{u_i} \biggr\}
  \prod_{i=1}^{n} (\pi x_i)^{2 \nu_n} 
  \prod_{i=1}^{n-1} u_i^{\nu_n} 
\\
& \qquad \cdot
  \widetilde{W}_{ \tilde{\nu} }^{n-1} 
  \Biggl(  x_2\sqrt{\frac{u_2}{u_1}}, \dotsc, 
          x_{n-1} \sqrt{\frac{u_{n-1}}{u_{n-2}}}, 
          x_n \sqrt{\frac{1}{u_{n-1}}} \Biggr)
  \prod_{i=1}^{n-1} \frac{du_i}{u_i}.
\end{split}
\end{equation}

The induction hypothesis implies the rapid decay of the function 
$ \widetilde{W}_{\nu}^{n-1}(y) $ and therefore the above integral converges absolutely 
for $ x \in (\R^+)^{n} $ and $ \nu \in \C^{n} $.
By using $ V_{\nu}(x) $, we can write $ \widetilde{W}_{\nu}^n(y) $ as
\begin{equation}
\begin{split}
  \widetilde{W}_{\nu}^n(y) 
&= \int_{(\R^+)^{n}}
    \prod_{i=1}^{n} \exp \biggl\{ -(\pi y_i)^2 t_i-\frac{1}{t_i} \biggr\}
    \prod_{i=1}^{n} t_i^{\nu_n}
\\
& \qquad \cdot
  V_{\nu} \Biggl( y_1 \sqrt{\frac{t_1}{t_2}}, \dotsc, 
                  y_{n-1} \sqrt{\frac{t_{n-1}}{t_n}}, 
                  y_n \sqrt{t_n} \Biggr) 
  \prod_{i=1}^{n} \frac{dt_i}{t_i}.
\end{split}
\end{equation}

We will first establish an expansion formula for
$ V_{\nu}(x) $ in Theorem \ref{expV} below. 
We notice that when $n=3$, this computation is essentially the same as \cite{HIO},
which we expressed  
the generalized principal series Whittaker functions on $Sp_3(\R)$
in terms of Whittaker functions on $SO_5(\R) $.
In the next step (subsection 3.2), by way of the results of Theorem \ref{expV} 
we will prove the relation (\ref{exp}).

\bigskip 
As in \cite{HIO}, to justify interchange of the order of 
integrations and infinite sums in the computation in the next subsections,
we need the following lemma. For the proof see \cite[Lemma 5.2]{HIO}.


\begin{lem} \label{conv}
For complex numbers $ \{ a_{ij} \}_{1 \le i\leq j \le n} $ 
with $ a_{ii} \geq 0  $, 
$ \{ b_i \}_{1\le i \le n} $ and $ d $, put 
$$
 \Delta(\bm) \equiv 
 \Delta(\bm, \{a_{ij} \}, \{ b_i \}, d)
 = \sum_{1\le i \le n} a_{ii} m_i^2 + \sum_{1 \le i < j \le n} a_{ij} m_i m_j 
   + \sum_{1 \le i \le n} b_i m_i + d
$$
Let $ \{ p_{ij} \}_{1 \le i,j \le n} $ be complex numbers and 
$ \{ q_i \}_{1 \le i \le n } $ nonzero complex numbers. 
We can inductively define complex numbers 
$ A_{\bm} \equiv A_{(m_1,\dotsc,m_n)}( \{a_{ij} \}, \{ b_i \}, d ) $ by 
$ A_{(0,\dotsc,0)} = 1 $ and 
the recurrence relation 
\begin{equation*} 
  \Delta(\bm) A_{\bm} 
 = \sum_{i=1}^n \Bigl( \sum_{j=1}^n  p_{ij} m_j + q_i \Bigr) A_{\bm-\be_i},
\end{equation*}
if $ \Delta(\bm) $ does not vanish for all $ (m_1,\dotsc,m_n) \neq (0,\dotsc,0) $.
Set 
$$ X= \{ ( \{a_{ij} \}, \{ b_i \}, d ) \in \C^{n(n+3)/2+1} \mid 
 \Delta(\bm) \neq 0 
  \mbox{ for all } \bm  \in \N^n \backslash \{(0,\dotsc,0)\}  \}.
$$ 
Let $U$ be any compact subset in $X$.
There exists a positive constant $ c_{U} $ depending only on $ U $ such that
\begin{equation} \label{esti}
 |A_{\bm}| \leq c_U^{m_1+ \cdots +m_n} / (m_1 + \cdots + m_n)! 
\end{equation}
for all $ \bm \in \N^n $  
and $ ( \{a_{ij} \}, \{ b_i \}, d ) \in U $.
Thus the power series 
$$\sum_{m_1,\dotsc,m_n=0}^{\infty} A_{\bm} x_1^{m_1} \cdots x_n^{m_n}$$
converges absolutely and uniformly on compacta for 
$ (x_1, \dotsc,x_n) \in (\R^+)^n $ and $ ( \{a_{ij} \}, \{ b_i \}, d ) \in X $. 
\end{lem}



\subsection{The first step --expansion formula for $V_{\nu} $--}

In this subsection we prove the following:
\begin{thm} \label{expV}
Let $ V_{\nu}(x) $ be the function defined by (\ref{defV}). 
Then, for $ \nu \in \,\!^{\prime} \!\!\> \a_{\C}^* $, we have
\begin{equation} \label{expformulaV}
  V_{\nu}(x) 
= \sum_{ w \in {\mathcal W}_{n-1} }
  \sum_{ 1 \le p \le n }
  w \Biggl[ \Gamma^{p} (\nu) 
  \sum_{\scriptstyle l_1, \dotsc, l_n=0 \atop \scriptstyle l_{n-1} \ge l_n}^{\infty} 
    b_{(l_1, \dotsc, l_n)}^{p} (\nu)
 \prod_{i=1}^n (\pi x_i)^{2(l_i+\nu_1^{(p)}+ \cdots + \nu_i^{(p)})} \Biggr],
\end{equation}
where 
\begin{align*}
 \Gamma^{p}(\nu)
& := \prod_{1 \le i<j \le n-1}  
     \Gamma ( -\nu_i-\nu_j )
     \Gamma ( -\nu_i+\nu_j )
   \prod_{1 \le i \le n-1} \Gamma(-2 \nu_i)
\\
& \quad \cdot
   \prod_{1 \le i \le p-1} \Gamma (-\nu_i+\nu_n)
   \prod_{p \le i \le n-1} \Gamma (\nu_i-\nu_n),
\end{align*}
$$ \nu^{(p)} \equiv (\nu_1^{(p)}, \dotsc, \nu_n^{(p)})
 := (\nu_1, \dotsc, \nu_{p-1},\nu_n, \nu_p, \dotsc,\nu_{n-1}). $$
and 
\begin{align*}
 b_{(l_1,\dotsc,l_n)}^{p}(\nu)
& := (-1)^{l_n} \sum_{ \{k_1, \dotsc, k_{n-2} \} }
   \frac{c_{n-1,(k_1,\dotsc,k_{n-2},l_n)}(\tilde{\nu})}
        { \prod_{  1 \le i \le \min ( p-1,n-2 )  }
              (l_i-k_i)! (\nu_i-\nu_n+1)_{l_i-k_{i-1}}  }
\\
& \qquad \qquad \cdot
   \frac{1}{ \prod_{ p \le i \le n-2  } 
              (l_i-k_{i-1})! (-\nu_i+\nu_n+1)_{l_i-k_i}  }
\\
& \qquad \qquad \cdot
\begin{cases}
  \frac{1}{(l_{n-1}-l_n)! (\nu_{n-1}-\nu_n+1)_{l_{n-1}-k_{n-2}} }
  & \mbox{ if } p = n,
\\
  \frac{1}{(l_{n-1}-k_{n-2})! (-\nu_{n-1}+\nu_n+1)_{l_{n-1}-l_n} }
  & \mbox{ if } p \neq n.
\end{cases}
\end{align*}
Here $ \{ k_1, \dotsc, k_{n-2} \} $ means $ k_i $ runs through such that
$$
    0 \le k_i \le l_i \ (1\le i \le p-2),
\ \ 0 \le k_{p-1} \le \min(l_{p-1},l_{p}), 
\ \ 0 \le k_i \le l_{i+1} \ (p \le i \le n-2).
$$
Moreover 
$ b_{\bl}^{p}(\nu) = b_{(l_1, \dotsc, l_n)}^{p} (\nu) $ 
is uniquely determined by the initial condition $ b_{(0,\dotsc,0)}^{p}(\nu) = 1 $ 
and the recurrence relation:
\[
 q_n(\bl, \nu^{(p)} ) b_{\bl}^{p}(\nu)
 = \sum_{i=1}^{n-1} b_{\bl-\be_i}^{p} (\nu)
+ \frac12 (-l_{n-1}+l_n+ \nu_{n}^{(p)}-\nu_n-1) b_{\bl-\be_n}^{p}(\nu).  
\]
\end{thm} 



\bpf
We substitute the expansion formula for $ W_{\tilde{\nu}}^{n-1}(y) $ to find 
\begin{align*}
 V_{\nu}(x) 
& = \sum_{w \in {\mathcal W}_{n-1}} 
    w \Biggl[
     \Gamma_{n-1}(\tilde{\nu})
  \int_{(\R^+)^{n-1}}
  \prod_{i=1}^{n-1} \exp \biggl\{ -(\pi x_i)^2 u_i - \frac{1}{u_i} \biggr\}
  \prod_{i=1}^{n} (\pi x_i)^{2\nu_n}
  \prod_{i=1}^{n-1} u_i^{\nu_n}
\\
& \qquad \cdot
  \sum_{k_1,\dotsc,k_{n-1}=0}^{\infty}
  c_{n-1, (k_1,\dotsc,k_{n-1})} (\tilde{\nu})
  \prod_{i=2}^{n} 
  \biggl(\pi x_i \sqrt{\frac{u_i}{u_{i-1}}}\biggr)^{\! 2(k_{i-1}+\nu_1+\cdots+\nu_{i-1})}
  \prod_{i=1}^{n-1} \frac{du_i}{u_i} \Biggr].
\end{align*}
By changing the order of the integration and the infinite sum, we get 
\begin{equation} \label{ZZ}
\begin{split}
 V_{\nu}(x) 
& = \sum_{w \in {\mathcal W}_{n-1}} 
    w \Biggl[
    \Gamma_{n-1}(\tilde{\nu})
   \sum_{k_1, \dotsc, k_{n-1}=0}^{\infty} 
    c_{n-1, (k_1,\dotsc,k_{n-1})}(\tilde{\nu})
   \prod_{i=1}^{n}
    (\pi x_i)^{2(k_{i-1}+\nu_1+ \cdots + \nu_{i-1} + \nu_n)}
\\
& \qquad \cdot
   \prod_{i=1}^{n-1}
   \int_{0}^{\infty}
       \exp \biggl\{ -(\pi x_i)^2 u_i - \frac{1}{u_i} \biggr\}
       u_i^{k_{i-1}-k_i- \nu_i +\nu_n}  \frac{du_i}{u_i}  
   \Biggr].
\end{split}
\end{equation}
As in \cite[\S 7]{HIO} this interchange is justified by Lemma \ref{conv},
and an analytic continuation argument implies that (\ref{ZZ}) is 
valid for all $ (x_1, \dotsc, x_n) \in (\R^+)^n $.

In view of 
\begin{equation} \label{Bessellem}
\begin{split}
 \int_0^{\infty}
 \exp \biggl\{ -(\pi x)^2 u - \frac{1}{u} \biggr\} u^{s} \frac{du}{u}
&= 2(\pi x)^{-s} K_{s}( 2\pi x)
\\
&= \frac{\pi}{\sin s \pi}
   (\pi x)^{-s} \bigl( I_{-s}(2\pi x) - I_{s}(2\pi x) \bigr)
\\
&= \frac{\pi}{\sin s \pi}
   \Biggl(
   \sum_{l =0}^{\infty} \frac{ (\pi x)^{2(l-s)} }{l! \,\Gamma(l-s+1) }
  -\sum_{l =0}^{\infty} \frac{ (\pi x)^{2l} }{ l! \,\Gamma(l+s+1) }
   \Biggr),
\end{split}
\end{equation}
we have
\begin{align*}
V_{\nu}(x) 
& = \sum_{w \in {\mathcal W}_{n-1}}
    w \Biggl[
    \Gamma_{n-1}(\tilde{\nu})
   \sum_{k_1, \dotsc, k_{n-1}=0 }^{\infty} 
    c_{n-1, (k_1,\dotsc,k_{n-1})}(\tilde{\nu})
   \prod_{i=1}^{n}
    (\pi x_i)^{2(k_{i-1}+\nu_1+ \cdots + \nu_{i-1} + \nu_n)}
\\
& \qquad \cdot
   \prod_{i=1}^{n-1}  \frac{\pi}{\sin (k_{i-1}-k_i-\nu_i+\nu_n) \pi }
\\
& \qquad \cdot
  \prod_{i=1}^{n-1}
  \Biggl(
   \sum_{l_i=0}^{\infty}
    \frac{ (\pi x_i)^{2(l_i-k_{i-1}+k_i + \nu_i -\nu_n) } }
         { l_i! \, \Gamma(l_i-k_{i-1}+k_i + \nu_i-\nu_n + 1) }
\\
& \qquad \qquad 
 - \sum_{l_i=0}^{\infty} 
    \frac{ (\pi x_i)^{2l_i} }
         { l_i! \, \Gamma(l_i+k_{i-1}-k_i - \nu_i+\nu_n + 1) }
  \Biggr) \Biggr]
\\
& = \sum_{w \in {\mathcal W}_{n-1}} 
    w \Biggl[
    \Gamma_{n-1}(\tilde{\nu})
   \sum_{k_1, \dotsc, k_{n-1}=0}^{\infty} 
    (-1)^{k_{n-1}}
    c_{n-1, (k_1,\dotsc,k_{n-1})}(\tilde{\nu})
    (\pi x_n)^{2(k_{n-1}+\nu_1+ \cdots + \nu_n)}
\\
& \qquad \cdot
  \prod_{i=1}^{n-1}
  \Biggl( 
   \Gamma( -\nu_i+\nu_n )
   \sum_{l_i=0}^{\infty} \frac{ (\pi x_i)^{2(l_i+k_i+\nu_1+\cdots+\nu_i)} }
                         { l_i! (\nu_i-\nu_n+1)_{l_i-k_{i-1}+k_i} }
\\
& \qquad \qquad 
+\Gamma( \nu_i-\nu_n )
   \sum_{l_i=0}^{\infty}
     \frac{ (\pi x_i)^{2(l_i+k_{i-1}+\nu_1+ \cdots + \nu_{i-1} + \nu_n)}}
                         {l_i! (-\nu_i+\nu_n+1)_{l_i+k_{i-1}-k_i} } 
   \Biggr) \Biggr].    
\end{align*}
By changing the order of the summuations and substituting 
$ l_i \to l_i-k_i $ or $ l_i \to l_i - k_{i-1} $ for $ 1 \le i \le n-1 $
and $ k_{n-1} \to l_n $ to get 
\begin{equation} \label{TTT}
\begin{split}
 V_{\nu}(x) 
&= \sum_{w \in {\mathcal W}_{n-1}}  \sum_{P \subset \{ 1,\dotsc, n-1 \} }
   w \Biggl[ 
   \Gamma^{P}(\nu) \sum_{l_1,\dotsc,l_n=0}^{\infty}
   b_{(l_1, \dotsc, l_n)}^{P}(\nu)
\\
& \qquad \cdot
  \prod_{i \in P \cup \{ n \} } (\pi x_i)^{2(l_i+\nu_1+\cdots+\nu_i)}
  \prod_{i \in P^c } (\pi x_i)^{2(l_i+\nu_1+\cdots+\nu_{i-1}+\nu_n)} \Biggr],
\end{split}
\end{equation}
where 
$ P $ ranges all the subset of $ \{1, \dotsc, n-1 \} $ and 
$ P^c $ means the complement of $ P $ in $ \{ 1, \dotsc, n-1 \} $.
Here 
\begin{align*}
 \Gamma^{P}(\nu)
:= \Gamma_{n-1}(\tilde{\nu})
  \prod_{i \in P}
     \Gamma ( -\nu_i+\nu_n)
  \prod_{i \in P^c}  
     \Gamma ( \nu_i-\nu_n)
\end{align*}
and
\begin{align*} 
 b_{(l_1,\dotsc,l_n)}^{P}(\nu)
& := (-1)^{l_n} \sum_{ \{k_1, \dotsc, k_{n-2} \} }
   \frac{c_{n-1, (k_1,\dotsc,k_{n-2},l_n)} (\tilde{\nu})}
        { \prod_{i \in P, 1 \le i \le n-2} 
              (l_i-k_i)! (\nu_i-\nu_n+1)_{l_i-k_{i-1}}  }
\\
& \qquad \cdot
   \frac{1}{ \prod_{i \in P^c, 1 \le i \le n-2} 
              (l_i-k_{i-1})! (-\nu_i+\nu_n+1)_{l_i-k_i}  }
\\
& \qquad \cdot
\begin{cases}
  \frac{1}{(l_{n-1}-l_n)! (\nu_{n-1}-\nu_n+1)_{l_{n-1}-k_{n-2}} }
  & \mbox{ if } n-1 \in P,
\\
  \frac{1}{(l_{n-1}-k_{n-2})! (-\nu_{n-1}+\nu_n+1)_{l_{n-1}-l_n} }
  & \mbox{ if } n-1 \in P^c,
\end{cases}
\end{align*}
for $ l_{n-1}-l_n \ge 0 $, and $ b_{(l_1,\dotsc,l_n)}^{P}(\nu) = 0 $ 
for $ l_{n-1}-l_n < 0 $.
Here the indexing set $ \{ k_1, \dotsc, k_{n-2} \} $ runs through such that
$$
  0 \le k_i \le l_i \ (i \in P),
\ \ 0 \le k_i \le l_{i+1} \ (i+1 \in P^c),
\ \ 0 \le k_i \ (i \in P^c \mbox{ and } i+1 \in P).
$$


From now on we consider which $ P $ contributes to the summation in (\ref{TTT}).
We first derive a recurrence relation for 
$ b_{(l_1,\dotsc,l_n)}^{P}(\nu) $ 
and an explicit formula for the initial value $ b_{(0,\dotsc,0)}^{P}(\nu) $.

\begin{lem} \label{lemmab}
For $ P \subset \{ 1,2,\dotsc,n-1 \} $, 
set $ \widetilde{P} 
 = \{ i \mid 1 \le i \le n-2, \ i \in P^c \mbox{ and } i+1 \in P \} $.
\begin{itemize}
\item[{\rm (i)}]
$ b_{\bl}^{P}(\nu) = b_{(l_1,\dotsc,l_n)}^{P}(\nu) $ 
satisfies the recurrence relation 
\begin{align*}
& \Bigl( \sum_{i=1}^{n-1} l_i^2 + \frac12 l_n^2 - \sum_{i=1}^{n-1} l_i l_{i+1}
  + \sum_{i=1}^{n-1} \lambda_i^P l_i + \kappa^P \Bigr)
 b_{\bl}^{P}(\nu)
\\
& = \sum_{i=1}^{n-1} b_{\bl-\be_i}^{P}(\nu)
   + \frac12 (-l_{n-1}+l_n+\lambda_n^P-\nu_n-1) b_{\bl-\be_n}^{P}(\nu)
\end{align*}
Here 
$$
  \kappa^{P} = \sum_{ i \in \widetilde{P} } (\nu_i-\nu_n)(\nu_{i+1}-\nu_n)
$$
and $ \lambda^P = (\lambda_{1}^{P},\dotsc,\lambda_n^P) $ is defined as follows:
\begin{align*}
\lambda_{i}^{P} = 
\begin{cases}
 \nu_i-\nu_{i+1}      & \mbox{ if } i-1 \in P,  i \in P,  i+1 \in P, \\
 \nu_i-\nu_n         & \mbox{ if } i-1 \in P,  i \in P,  i+1 \in P^c, \\
 \nu_{i-1}+\nu_i-\nu_{i+1}-\nu_n 
                          & \mbox{ if } i-1 \in P^c, i \in P, i+1 \in P, \\
 \nu_{i-1}+\nu_i-2\nu_n & \mbox{ if } i-1 \in P^c,  i \in P,i+1 \in P^c, \\
 -\nu_i-\nu_{i+1}+2\nu_n & \mbox{ if } i-1 \in P,  i \in P^c,  i+1 \in P, \\
 -\nu_i+\nu_n  & \mbox{ if } i-1 \in P,  i \in P^c,  i+1 \in P^c, \\
 \nu_{i-1}-\nu_i-\nu_{i+1}+\nu_n
               & \mbox{ if } i-1 \in P^c,  i \in P^c,  i+1 \in P,\\
 \nu_{i-1}-\nu_i & \mbox{ if } i-1 \in P^c,  i \in P^c,  i+1 \in P^c, 
\end{cases} 
\end{align*}
for $ 1 \le i \le n-2 $,
\begin{align*}
 \lambda_{n-1}^{P} =
\begin{cases}
 \nu_{n-1}-\nu_n            & \mbox{ if } n-2 \in P, \ n-1 \in P, \\
 \nu_{n-2}+\nu_{n-1}-2\nu_n & \mbox{ if } n-2 \in P^c, \ n-1 \in P, \\
 -\nu_{n-1}+\nu_n           & \mbox{ if } n-2 \in P, \ n-1 \in P^c, \\
 \nu_{n-2}-\nu_{n-1}        & \mbox{ if } n-2 \in P^c, \ n-1 \in P^c,
\end{cases}
\end{align*}
and
\begin{align*} 
  \lambda_n^{P}  =
\begin{cases} 
   \nu_n     & \mbox{ if } n-1 \in P,  \\ 
   \nu_{n-1} & \mbox{ if } n-1 \in P^c.
\end{cases}
\end{align*} 
\item[{\rm (ii)}]
We have 
\begin{align*}
 b_{(0,\dotsc,0)}^{P}(\nu)
= \prod_{i \in \widetilde{P} }
   \frac{ \Gamma( \nu_i-\nu_{i+1}+ 1 ) }
        { \Gamma( \nu_i-\nu_n + 1 ) \Gamma( -\nu_{i+1}+\nu_n + 1) }
\end{align*}
and $ b_{(0,\dotsc,0)}^{P}(\nu) = 1 $ when 
$ \widetilde{P} = \emptyset $.
\end{itemize}
\end{lem} 
\bpf
(i) The idea of proof is similar to the first step in 
the proof of Theorem \ref{fund} $(P=\{ 1, \dotsc, n-1 \})$.
Our claim follows from the identity
\begin{align*}
& \sum_{i \in P, \, 1 \le i \le n-2}
   (l_i-k_i) (l_i-k_{i-1}+ \nu_i-\nu_n)
 +  \sum_{i \in P^c, \, 1 \le i \le n-2 }
   (l_i-k_{i-1}) (l_i-k_i-\nu_i+\nu_n)
\\
& \quad + \begin{cases}
   (l_{n-1}-l_n)(l_{n-1}-k_{n-2}+\nu_{n-1}-\nu_n) & \mbox{ if } n-1 \in P, \\
   (l_{n-1}-k_{n-2})(l_{n-1}-l_n-\nu_{n-1}+\nu_n) & \mbox{ if } n-1 \in P^c 
   \end{cases}
\\
& \quad - \sum_{i \in P, i+1 \in P, \, 1 \le i \le n-2}
         (l_i-k_i) (l_{i+1}-k_i+ \nu_{i+1}-\nu_n)
\\
& \quad
   - \sum_{i \in P, i+1 \in P^c, \, 1 \le i \le n-2}
         (l_i-k_i)(l_{i+1}-k_i)
\\
& \quad - \sum_{i \in P^c, i+1 \in P, \, 1 \le i \le n-2}
       (l_i-k_i-\nu_i+\nu_n) 
       (l_{i+1}-k_i+\nu_{i+1}-\nu_n)
\\
& \quad
   - \sum_{i \in P^c, i+1 \in P^c, \, 1 \le i \le n-2}
       (l_i-k_i-\nu_i+\nu_n) (l_{i+1}-k_i)
\\
&\quad = \Bigl( \sum_{i=1}^{n-1} l_i^2 + \frac12 l_n^2 - \sum_{i=1}^{n-1} l_i l_{i+1}
  + \sum_{i=1}^{n-1} \lambda_i^P l_i + \kappa^P \Bigr)
    -q_{n-1} \bigl( (k_1,\dotsc,k_{n-2},l_n), \tilde{\nu} \bigr).
\end{align*}

\noindent
(ii) From the definition of $ b_{(l_1,\dotsc,l_n)}^{P}(\nu) $, 
\begin{align*}
b_{(0, \dotsc, 0)}^{P}(\nu)
& = \sum_{ \scriptstyle k_i=0 \atop \scriptstyle (i \in \widetilde{P}) }^{\infty}
   \frac{ c_{n-1, (k_1,\dotsc, k_{n-2},0)}^{\widetilde{P}}(\nu) }
        { \prod_{i \in P \cap \{1,\dotsc,n-2\} } 
          (\nu_i-\nu_n+1)_{-k_{i-1}}  
          \prod_{i \in P^c \cap \{1,\dotsc,n-2\} } 
          (-\nu_i+\nu_n+1)_{-k_i} } 
\\ & \qquad \qquad 
   \cdot 
   \begin{cases}
   \frac{1}{ (\nu_{n-1}-\nu_n+1)_{-k_{n-2}} } 
         & \mbox{ if } n-2 \in \widetilde{P}, \\
   1 & \mbox{ otherwise }
   \end{cases}
\\
& = \sum_{ \scriptstyle k_i=0 \atop \scriptstyle (i \in \widetilde{P}) }^{\infty}
    \Biggl( \prod_{i \in \widetilde{P} }
    \frac{1}{ (\nu_{i+1}-\nu_n+1)_{-k_{i}} (-\nu_i+\nu_n+1)_{-k_i} }
    \Biggr) 
    \cdot c_{n-1,(k_1,\dotsc, k_{n-2},0)}^{\widetilde{P}}(\nu),
\end{align*}
where
$$  c_{n-1, (k_1,\dotsc, k_{n-2},0)}^{\widetilde{P}}(\nu) = 
   c_{n-1, (k_1,\dotsc, k_{n-2},0)}(\nu)|_{k_i = 0 \ (i \notin \widetilde{P})}. $$
Since   
$ i \in \widetilde{P} $ implies $  i \pm 1 \notin \widetilde{P} $,
we can find the recurrence relation for 
$  c_{n-1,(k_1,\dotsc, k_{n-2},0)}^{\widetilde{P}}(\nu) $:
\begin{align*}
 \biggl( \sum_{i \in \widetilde{P} } k_i^2 
 +\sum_{i \in \widetilde{P}} (\nu_i-\nu_{i+1}) k_i \biggr) 
   c_{n-1, (k_1,\dotsc, k_{n-2},0)}^{\widetilde{P}}(\nu)
 = \sum_{i \in \widetilde{P}} 
   c_{n-1, (k_1,\dotsc, k_{n-2},0)-{\be_i} }^{\widetilde{P}}(\nu).
\end{align*}
We can easily solve it to find
$$
  c_{n-1, (k_1,\dotsc, k_{n-2},0)}^{\widetilde{P}}(\nu)
 = \prod_{i \in \widetilde{P} } \frac{1}{k_i! (\nu_i-\nu_{i+1}+1)_{k_i} }.
$$
Thus we get
\begin{align*}
 b_{(0, \dotsc, 0)}^{P}(\nu)
& = \prod_{i\in \widetilde{P}} \biggl( \sum_{k_i=0}^{\infty}
    \frac{1}{ k_i! (\nu_i-\nu_{i+1}+1)_{k_i} 
             (\nu_{i+1}-\nu_n+1)_{-k_{i}} (-\nu_i+\nu_n+1)_{-k_i} }
     \biggr)
\\
& = \prod_{i \in \widetilde{P} }
    \,_2F_1 ( -\nu_{i+1}+\nu_n, \nu_i-\nu_n ;  \nu_i-\nu_{i+1} ; 1 ) 
\\
& = \prod_{i \in \widetilde{P} }
    \frac{ \Gamma( \nu_i-\nu_{i+1} + 1 ) }
        { \Gamma( \nu_i-\nu_n + 1 ) \Gamma( -\nu_{i+1}+\nu_n + 1) }.
\end{align*}
Here we used $ (a)_{-n} = (-1)^n / (1-a)_n $ and Gauss' formula
$$
 _2F_1(a,b;c;1)
 = \frac{ \Gamma(c-a-b) \Gamma(c) }{ \Gamma(c-a) \Gamma(c-b) }
$$
in the last step.
Therefore we complete the proof of Lemma \ref{lemmab}. 
\epf

\bigskip

Returning to the proof of Theorem \ref{expV},
the next proposition implies our expansion formula 
for $ V_{\nu}(x) $.
Because, 
for $ 1\le p \le n $, we have
$$ 
 \Gamma^{\{1,2, \dotsc, p-1 \}}(\nu) = \Gamma^p(\nu), \ \ 
 b_{(l_1,\dotsc,l_n)}^{ \{1,2, \dotsc, p-1 \}}(\nu) = b_{(l_1,\dotsc,l_n)}^p(\nu),
$$
from the definition.

\begin{prop} \label{claimV}
The following $P$ contributes in the right hand side of (\ref{TTT}).
\begin{itemize}
\item $ P = \emptyset $,
\item $ P $ is of the form $ \{1,2, \dotsc, p-1 \} $ for some $ 2 \le p \le n $.
\end{itemize}
More precisely, 
if there exists an element $ p_0 $ in $ \widetilde{P} $ then we have
\begin{align*}
& \sum_{ w \in \{ 1, w_{p_0}  \} } 
 w \Biggl[ \Gamma^P(\nu)
 \sum_{l_1, \dotsc, l_n =0}^{\infty} 
 b_{(l_1,\dotsc, l_n)}^{P}(\nu) 
\\
& \quad \cdot
 \prod_{i \in P \cup \{n \} } (\pi x_i)^{2(l_i+\nu_1+\cdots+\nu_i)} 
 \prod_{i \in P^c } (\pi x_i)^{2(l_i+\nu_1+\cdots+\nu_{i-1}+\nu_n)} \Biggr]
 = 0.
\end{align*}
Here $ w_{p_0} $ is the simple reflection of the Weyl group $ {\mathcal W}_{n-1} $,
that is, it
permutes $ \nu_{p_0} $ and $ \nu_{p_0+1} $.
\end{prop}


\bpf
Fix $ p_0 \in \widetilde{P} $.
Since $ p_0 \in P^c $ and $ p_0+1 \in P $, 
$$
 \prod_{i \in P \cup \{n \} } (\pi x_i)^{2(l_i+\nu_1+\cdots+\nu_i)} 
 \prod_{i \in P^c } (\pi x_i)^{2(l_i+\nu_1+\cdots+\nu_{i-1}+\nu_n)} 
$$
is invariant under the permutation of $ \nu_{p_0} $ and $ \nu_{p_0+1} $.
Then it is enough to show 
$$
 a_{(l_1,\dotsc,l_n)}^{P}(\nu)
 := \sum_{ w \in \{ 1, w_{p_0} \} }
  w \Bigl[ \Gamma^{P}(\nu) b_{(l_1,\dotsc, l_n)}^{P}(\nu) \Bigr]
  = 0.
$$
In view of Lemma \ref{lemmab} (i), 
we can check that $ \lambda_i^{P} $ and $ \kappa^{P} $
is invariant under the action of $ w_{p_0} $.
Then 
$ b_{(l_1,\dotsc, l_n)}^{P}(\nu) $ and 
$ b_{(l_1,\dotsc, l_n)}^{P}(w_{p_0} \nu) $ 
satisfies the same recurrence relation and therefore 
$ a_{(l_1,\dotsc,l_n)}^{P}(\nu) $ also satisfies the same one.
Thus, if we can say 
$ a_{(0,\dotsc,0)}^{P}(\nu) = 0 $ then 
$ a_{(l_1,\dotsc,l_n)}^{P}(\nu) = 0$ inductively follows.

From Lemma \ref{lemmab} (ii), we have 
\begin{align*}
\Gamma^p(\nu) b_{(0,\dotsc,0)}^{P}(\nu)
& = 
  \prod_{1 \le i<j \le n-1} 
     \Gamma ( -\nu_i-\nu_j)
     \Gamma ( -\nu_i+\nu_j)
  \prod_{1 \le i \le n-1} \Gamma(-2\nu_i) 
\\
& \quad \cdot 
  \prod_{i \in P}
     \Gamma ( -\nu_i+\nu_n )
  \prod_{i \in P^c}  
     \Gamma ( \nu_i-\nu_n )
  \prod_{i \in \widetilde{P} }
   \frac{ \Gamma( \nu_i-\nu_{i+1} + 1 ) }
        { \Gamma( \nu_i-\nu_n + 1 ) \Gamma( -\nu_{i+1}+\nu_n + 1) }.
\end{align*}
We pick up the terms which are not invariant under the action of 
$ w_{p_0} $:
\begin{align*}
& \Gamma ( -\nu_{p_0}+\nu_{p_0+1} ) 
  \Gamma ( -\nu_{p_0+1} -\nu_n )
  \Gamma ( \nu_{p_0}-\nu_n )
  \cdot \frac{ \Gamma(\nu_{p_0}-\nu_{p_0+1}+1)}
      {\Gamma( \nu_{p_0}-\nu_n+1) \Gamma(-\nu_{p_0+1}+\nu_n+1) }
\\
&= \frac{\pi}{\sin (-\nu_{p_0}+\nu_{p_0+1}) \pi }
   \cdot \frac{1}{ (\nu_{p_0}-\nu_n)( -\nu_{p_0+1}+\nu_n) }.
\end{align*}
Therefore we get 
$ a_{(0,\dotsc, 0)}^{P}(\nu) = 0$
and complete the proof of Theorem \ref{expV}.
\epf

\subsection{The second step --expansion formula for $\widetilde{W}_{\nu}^n$--}
In the similar way to the previous subsection, we shall prove the linear 
relation (\ref{exp}). We need a little more complicated argument.
We insert the expansion formula (\ref{expformulaV})
for $ V_{\nu}(x) $ to get
\begin{align*}
\widetilde{W}_{\nu}^n(y)
& =
\sum_{ w \in {\mathcal W}_{n-1} } \sum_{p=1}^{n}
  w \Biggl[ 
  \Gamma^{p}(\nu) \int_{(\R^+)^n}
  \prod_{i=1}^{n} \exp \biggl\{ -(\pi y_i)^2 t_i - \frac{1}{t_i} \biggr\}
  \prod_{i=1}^{n} t_i^{\nu_n}
\\
& \qquad \cdot
  \sum_{\scriptstyle l_1, \dotsc, l_n=0 \atop \scriptstyle l_{n-1} \ge l_n}^{\infty} 
 b_{(l_1,\dotsc,l_n)}^{p}(\nu)
  \prod_{i=1}^{n} \biggl( \pi y_i \sqrt{\frac{t_i}{t_{i+1}}}\biggr)^{\! 
                  2(l_i + \nu_1^{(p)} + \cdots + \nu_i^{(p)}) } 
  \prod_{i=1}^{n} \frac{dt_i}{t_i} \Biggr].
\end{align*}
By changing the order of the integration and the infinite sum, we have
\begin{align*}
\widetilde{W}_{\nu}^n(y)
& = 
\sum_{ w \in {\mathcal W}_{n-1} } \sum_{p=1}^{n}
  w \Biggl[ 
  \Gamma^{p}(\nu) 
  \sum_{\scriptstyle l_1, \dotsc, l_n=0 \atop \scriptstyle l_{n-1} \ge l_n}^{\infty} 
   b_{(l_1,\dotsc,l_n)}^{p}(\nu)
  \prod_{i=1}^{n} ( \pi y_i )^{2(l_i + \nu_1^{(p)} + \cdots + \nu_i^{(p)}) }
\\
& \qquad \cdot
  \prod_{i=1}^{n}
  \int_{0}^{\infty}
  \exp \biggl\{ -(\pi y_i)^2 t_i-\frac{1}{t_i} \biggr\} 
 \, t_i^{-l_{i-1}+l_i+\nu_i^{(p)}+\nu_n} \frac{dt_i}{t_i} \Biggr].
\end{align*}
We use (\ref{Bessellem}) for the integral above to find
\begin{align*}
\widetilde{W}_{\nu}^n(y)
& = 
\sum_{ w \in {\mathcal W}_{n-1} } \sum_{p=1}^{n}
  w \Biggl[ 
  \Gamma^{p}(\nu) 
  \sum_{\scriptstyle l_1, \dotsc, l_n=0 \atop \scriptstyle l_{n-1} \ge l_n}^{\infty} 
  (-1)^{l_n} b_{(l_1,\dotsc,l_n)}^{p}(\nu)
\\
& \quad \cdot 
  \prod_{i=1}^{n} 
  \Biggl\{
   \Gamma (- \nu_i^{(p)} - \nu_n )
   \sum_{m_i=0}^{\infty}
   \frac{ ( \pi y_i )^{ 2(m_i+l_i + 
           \nu_1^{(p)} + \cdots + \nu_{i}^{(p)}) }}
      {m_i! ( \nu_i^{(p)}+\nu_n + 1 )_{m_i-l_{i-1}+l_i} } 
\\
& \qquad 
 + \Gamma ( \nu_i^{(p)} + \nu_n )
   \sum_{m_i=0}^{\infty}
   \frac{ ( \pi y_i )^{ 2(m_i+l_{i-1} + 
           \nu_1^{(p)} + \cdots + \nu_{i-1}^{(p)} - \nu_n)}}
        {m_i! (- \nu_i^{(p)} -\nu_n + 1 )_{m_i+l_{i-1}-l_i} }
 \Biggr\} \Biggr].
\end{align*}

\noindent
We substitute $ m_i \to m_i - l_{i-1} $ or $ m_i \to m_i-l_i $, 
and arrange the order of the summation. Then we get
\begin{equation} \label{UUU}
\begin{split}
\widetilde{W}_{\nu}^n(y)
&= 
\sum_{ w \in {\mathcal W}_{n-1} } \sum_{p=1}^{n} \sum_{Q \subset \{1,2,\dotsc,n \} }
w \Bigl[ 
  \Gamma^{p,Q}(\nu) 
  \widetilde{M}_{\nu}^{p,Q}(y)
  \Bigr],
\end{split}
\end{equation}
where
\begin{align*}
 \Gamma^{p,Q}(\nu)
& := \Gamma^{p}(\nu)
    \prod_{i \in Q} \Gamma ( -\nu_i^{(p)} - \nu_n)
    \prod_{i \in Q^c} \Gamma (  \nu_i^{(p)} + \nu_n )
\end{align*}
and
\begin{align*}
 \widetilde{M}_{\nu}^{p,Q}(y)
 := \sum_{m_1, \dotsc, m_n=0}^{\infty}  c_{(m_1,\dotsc,m_n)}^{p,Q}(\nu)
  \prod_{i \in Q} (\pi y_i)^{2(m_i+\nu_1^{(p)} + \cdots + \nu_{i}^{(p)} )}
  \prod_{i \in Q^c} (\pi y_i)^{2(m_i+\nu_1^{(p)} + \cdots + \nu_{i-1}^{(p)} -\nu_n) } 
\end{align*}
is the power series with the coefficient  
\begin{align*}
 c_{(m_1,\dotsc,m_n)}^{p,Q}(\nu)
& := \sum_{ \{ l_1, \dotsc, l_{n} \} } (-1)^{l_n} b_{(l_1,\dotsc,l_n)}^{p} (\nu)
   \prod_{i \in Q}
  \frac{1}{(m_i-l_i)! ( \nu_i^{(p)}+\nu_n +1 )_{m_i-l_{i-1}}} 
\\
& \qquad \cdot
  \prod_{i \in Q^c} 
  \frac{1}{(m_i-l_{i-1})! ( -\nu_i^{(p)}-\nu_n +1 )_{m_i-l_i}}.
\end{align*}
Here $ \{ l_1, \dotsc, l_n \} $ means that
$$
  0 \le l_i \le m_i \ ( i \in Q ),
\ \ 0 \le l_i \le m_{i+1} \ (i+1 \in Q^c), 
\ \ 0 \le l_i \ ( i \in Q^c \mbox{ and } i+1 \in Q).
$$

As in the previous subsection
let us derive a recurrence relation
for $ c^{p,Q}_{(m_1,\dotsc,m_n)}(\nu) $ and an explicit formula for 
the initial value $ c_{(0,\dotsc,0)}^{p,Q}(\nu) $.

\begin{lem} \label{recU}
For $ Q \subset \{ 1,2, \dotsc, n \} $, set 
$ \widetilde{Q} = \{ i \mid 1 \le i \le n-1, \ i \in Q^c \mbox{ and } i+1 \in Q \} $.
\begin{itemize}
\item[{\rm (i)}]
$ c_{\bm}^{p,Q}(\nu) = c_{(m_1,\dotsc,m_n)}^{p,Q}(\nu) $ satisfies the 
recurrence relation
\begin{align*}
 \Bigl( \sum_{i=1}^{n-1} m_i^2 + \frac12 m_n^2 - \sum_{i=1}^{n-1} m_i m_{i+1}
  + \sum_{i=1}^{n-1} \lambda_{i}^{p,Q} m_i + \kappa^{p,Q} \Bigr)
  c_{\bm}^{p,Q}(\nu) 
= \sum_{i=1}^{n-1} c_{\bm-\be_i}^{p,Q}(\nu) + \frac12 c_{\bm-\be_n}^{p,Q}(\nu),
\end{align*}
where
\begin{equation*}
 \kappa^{p,Q} 
= \sum_{i \in \widetilde{Q} }
  (\nu_i^{(p)} + \nu_n)( \nu_{i+1}^{(p)} + \nu_n ) 
 + \begin{cases}
    0 & \mbox{ if } n \in Q, \\
    \frac12 ( \nu_{n}^{(p)} + \nu_n )( -\nu_{n}^{(p)} + \nu_n ) & \mbox{ if } n \in Q^c,
   \end{cases}
\end{equation*}
and $ \lambda^{p,Q} = (\lambda_1^{p,Q}, \dotsc, \lambda_n^{p,Q}) $ is defined as follows:
\begin{equation*}
 \lambda_i^{p,Q} =
\begin{cases}
 \nu_{i}^{(p)} - \nu_{i+1}^{(p)}
   & \mbox{ if } i-1 \in Q, \ i \in Q, \ i+1 \in  Q, \\
 \nu_{i}^{(p)} + \nu_n
   & \mbox{ if } i-1 \in Q, \ i \in Q, \ i+1 \in  Q^c, \\
 \nu_{i-1}^{(p)} + \nu_{i}^{(p)} - \nu_{i+1}^{(p)} + \nu_n
   & \mbox{ if } i-1 \in Q^c, \ i \in Q, \ i+1 \in  Q, \\
 \nu_{i-1}^{(p)} + \nu_{i}^{(p)} + 2 \nu_n
   & \mbox{ if } i-1 \in Q^c, \ i \in Q, \ i+1 \in  Q^c, \\
 -\nu_{i}^{(p)} - \nu_{i+1}^{(p)} - 2\nu_n
   & \mbox{ if } i-1 \in Q, \ i \in Q^c, \ i+1 \in  Q, \\
 -\nu_{i}^{(p)} - \nu_n
   & \mbox{ if } i-1 \in Q, \ i \in Q^c, \ i+1 \in  Q^c, \\
 \nu_{i-1}^{(p)} - \nu_{i}^{(p)} - \nu_{i+1}^{(p)} - \nu_n
   & \mbox{ if } i-1 \in Q^c, \ i \in Q^c, \ i+1 \in  Q, \\
 \nu_{i-1}^{(p)} - \nu_{i}^{(p)}  
   & \mbox{ if } i-1 \in Q^c, \ i \in Q^c, \ i+1 \in  Q^c,
\end{cases}
\end{equation*}
for $ 1 \le i \le n-1 $ and 
\begin{equation*}
 \lambda_n^{p,Q} = 
\begin{cases}
  \nu_{n}^{(p)}    & \mbox{ if } n-1 \in Q, \ n \in Q, \\ 
  \nu_{n-1}^{(p)} + \nu_{n}^{(p)} + \nu_n  & \mbox{ if } n-1 \in Q^c, \ n \in Q, \\
  -\nu_n
                           & \mbox{ if } n-1 \in Q, \ n \in Q^c, \\
  \nu_{n-1}^{(p)}     & \mbox{ if } n-1 \in Q^c, \ n \in Q^c.
\end{cases}
\end{equation*}
\item[{\rm (ii)}]
We have 
\begin{align*}
 c_{(0,\dotsc,0)}^{p,Q}(\nu)
& = \prod_{i \in \widetilde{Q} } 
   \frac{\Gamma ( \nu_i^{(p)} - \nu_{i+1}^{(p)} + 1 ) }
        {\Gamma ( -\nu_{i+1}^{(p)} - \nu_n + 1 ) 
         \Gamma ( \nu_i^{(p)} + \nu_n + 1 )}
\\
& \qquad \cdot 
  \begin{cases}
    1 & \mbox{ if } n \in Q, \\ 
    \dfrac{\Gamma(2 \nu_n^{(p)}+1)}
         {\Gamma(\nu_n^{(p)}+\nu_n+1 )
          \Gamma(\nu_n^{(p)}-\nu_n+1 ) } &  \mbox{ if } n \in Q^c.
  \end{cases}
\end{align*}
\end{itemize}
\end{lem}
\bpf
(i) 
Our claim follows from the identity
\begin{align*}
& \sum_{i \in Q, \, 1 \le i \le n-1 } 
  (m_i-l_{i}) (m_i-l_{i-1} + \nu_{i}^{(p)}+\nu_n)
+\sum_{i \in Q^c, \, 1 \le i \le n-1 } 
  (m_i-l_{i-1}) (m_i-l_i- \nu_{i}^{(p)}-\nu_n)
\\
& \ \ + \frac12
  \begin{cases}
  (m_n-l_n) (m_n-l_{n-1} + \nu_{n}^{(p)} +\nu_n) 
       &  \mbox{ if } n \in Q, \\
  (m_n - l_{n-1}) (m_n-l_n -\nu_{n}^{(p)}-\nu_n) 
       &  \mbox{ if } n \in Q^c
  \end{cases}
\\
& \ \  
- \sum_{i \in Q, i+1 \in Q, \, 1 \le i \le n-1} 
  (m_i-l_i) (m_{i+1}-l_i + \nu_{i+1}^{(p)}+\nu_n)
\\
& \ \
- \sum_{i \in Q, i+1 \in Q^c, \, 1 \le i \le n-1} 
  (m_i-l_i)(m_{i+1}-l_i)
\\
& \ \
- \sum_{i \in Q^c, i+1 \in Q,  \, 1 \le i \le n-1} 
  (m_i-l_{i} - \nu_{i}^{(p)}-\nu_n )
  (m_{i+1}-l_i + \nu_{i+1}^{(p)}+\nu_n )
\\
& \ \ 
- \sum_{i \in Q^c, i+1 \in Q^c, \, 1 \le i \le n-1} 
  (m_i-l_{i} - \nu_{i}^{(p)}-\nu_n)(m_{i+1}-l_i)
\\
& \ \  -\frac12
  \begin{cases}
   (m_n-l_n - \nu_{n}^{(p)} - \nu_n) 
  (l_{n-1}-l_n - \nu_{n}^{(p)} + \nu_n) 
     & \mbox{ if } n \in Q, \\
  (m_n-l_n) (l_{n-1}-l_n - \nu_{n}^{(p)}+\nu_n)
     & \mbox{ if } n \in Q^c \\
  \end{cases}
\\
& \ \ = \Bigl( \sum_{i=1}^{n-1} m_i^2 + \frac12 m_n^2 - \sum_{i=1}^{n-1} m_i m_{i+1}
  + \sum_{i=1}^{n-1} \lambda_{i}^{p,Q} m_i + \kappa^{p,Q} \Bigr)
  - q_n(\bl,\nu^{(p)}).
\end{align*}
(ii) 
We can prove in the same way as Lemma \ref{lemmab} (ii).
\epf

\bigskip 

The following is immediate from the above lemma.

\begin{cor} \label{cor1}
We have 
\begin{itemize}
\item 
$ \widetilde{M}_{\nu}^{p, \{1,2,\dotsc,n\}}(y) = \widetilde{M}_{\nu^{(p)}}^n(y) $ and
$ \Gamma^{p,\{1,2, \dotsc,n \}}(\nu) = \Gamma_{n}(\nu^{(p)}) $ 
for $  1 \le p \le n $,
\item 
$ \widetilde{M}_{\nu}^{n, \{1,2,\dotsc,q-1\}}(y) = \widetilde{M}_{\bar{\nu}^{(q)}}(y) $ and
$ \Gamma^{n,\{1,2,\dotsc,q-1\}}(\nu) = \Gamma_{n}(\bar{\nu}^{(q)}) $ 
for $ 1 \le q \le n $. 
\end{itemize}
Here we write 
$ \bar{\nu}^{(q)} := (\nu_1,\dotsc,\nu_{q-1},-\nu_n,\nu_{q},\dotsc,\nu_{n-1}) $.
\end{cor}

\bpf 
By Lemma \ref{recU} (i), we can verify
$ c_{\bm}^{p, \{1,2, \dotsc,n \}}(\nu) $ $ ( 1 \le p \le n) $ and 
$ c_{\bm}^{n,\{1,2,\dotsc,q-1 \}}(\nu) $ $ (1 \le q \le n) $ 
satisfy the same recurrence
relations as (\ref{rec}) with $ \nu = \nu^{(p)} $
and 
$ \nu =\bar{\nu}^{(q)} $, respectively.
Since
$ c_{(0,\dotsc,0)}^{p, \{1,2, \dotsc,n \}}(\nu) 
 = c_{(0,\dotsc,0)}^{n,\{1,2, \dotsc,q-1\}} (\nu) = 1 $ 
from Lemma \ref{recU} (ii), we have
$ c_{\bm}^{p, \{1,2, \dotsc,n \} }(\nu) = c_{n,\bm}(\nu^{(p)}) $ 
and $ c_{\bm}^{n,\{1,2,\dotsc,q-1\}}(\nu) = c_{n,\bm}( \bar{\nu}^{(q)}  ) $ 
and thus get our claim for $ \widetilde{M}^{*,*}_{\nu}(y) $.

The latter can be seen from the definition. 
Indeed we have
\begin{align*}
 \Gamma^{p,\{1,2, \dotsc,n \}}(\nu) 
& = \Gamma_{n-1}(\tilde{\nu}) 
    \prod_{1 \le i \le p-1} \Gamma(-\nu_i+\nu_n) 
    \prod_{p \le i \le n-1} \Gamma(\nu_i-\nu_n) 
    \prod_{1 \le i \le n} \Gamma(-\nu_i^{(p)}-\nu_n) 
\\
& =\prod_{1 \le i<j \le n-1} 
      \Gamma(-\nu_i-\nu_j) \Gamma(-\nu_i+\nu_j)
   \prod_{1 \le i \le n-1} \Gamma(-2\nu_i) 
\\
& \quad \cdot 
   \prod_{1 \le i \le p-1} \Gamma(-\nu_i+\nu_n) 
   \prod_{p \le i \le n-1} \Gamma(\nu_i-\nu_n) 
\\
& \quad \cdot 
    \prod_{1 \le i \le p-1} \Gamma(-\nu_i-\nu_n) \cdot \Gamma(-\nu_n-\nu_n) \cdot 
    \prod_{p+1 \le i \le n} \Gamma(-\nu_{i-1}-\nu_n)
\\
& = \prod_{1 \le i<j \le n-1} 
     \Gamma(-\nu_i-\nu_j) \Gamma(-\nu_i+\nu_j)
    \prod_{1 \le i \le n} \Gamma(-2\nu_i) 
\\
& \quad \cdot 
    \prod_{1 \le i \le p-1} \Gamma(-\nu_i+\nu_n) \Gamma(-\nu_i-\nu_n)
    \prod_{p+1 \le i \le n} \Gamma(\nu_{i-1}-\nu_n) \Gamma(-\nu_{i-1}-\nu_n)
\\
& = \Gamma_{n}(\nu^{(p)}) 
\end{align*} 
and
\begin{align*}
& \Gamma^{n,\{1,2,\dotsc,q-1\}}(\nu) 
\\
& = \Gamma_{n-1}(\tilde{\nu})
    \prod_{1 \le i \le n-1} \Gamma(-\nu_i+\nu_n)   
    \prod_{1 \le i \le q-1} \Gamma(-\nu_i^{(n)}-\nu_n) 
    \prod_{q \le i \le n} \Gamma(\nu_i^{(n)}+\nu_n)  
\\
& = \prod_{1 \le i \le n-1} \Gamma(-2\nu_i) \cdot \Gamma(2\nu_n) 
  \prod_{1 \le i<j \le n-1} 
    \Gamma(-\nu_i-\nu_j) \Gamma(-\nu_i+\nu_j)
\\
& \quad \cdot 
  \prod_{1 \le i \le q-1} \Gamma(-\nu_i+\nu_n) \Gamma(-\nu_i-\nu_n) 
  \prod_{q+1 \le i \le n} \Gamma(-\nu_{i-1}+\nu_n)\Gamma(\nu_{i-1}+\nu_n)
\\
& = \Gamma_{n}(\bar{\nu}^{(q)}) .
\end{align*}
\epf

\bigskip 

Now we state cancelations in the summation in (\ref{UUU}).


\begin{prop} \label{claimU}
The following terms contribute to the right hand side of (\ref{UUU}).
\begin{itemize}
\item $ 1 \le p \le n $ and $ Q= \{1,2,\dotsc,n\}  $,
\item $ p=n$ and $ Q $ is of the form $ \{ 1,2, \dotsc,q-1 \} $ $ (1 \le q \le n) $.
\end{itemize}
More precisely we have the following. 
\begin{itemize}
\item[{\rm (i)}]
If $ 1 \le p \le n-1 $
and $  Q $ is of the form $ \{ 1,2, \dotsc,q-1 \} $ $ (1 \le q \le n) $, 
then   
\begin{align*}
& \sum_{w \in \{ 1, w_{n-1} \} }
  w \Bigl[ \Gamma^{p,Q}(\nu) \widetilde{M}_{\nu}^{p,Q}(y) \Bigr] = 0.
\end{align*}
Here $ w_{n-1} \in {\mathcal W}_{n-1} $ is the simple reflection, which
permutes the sign of $ \nu_{n-1} $.

\item[{\rm (ii)}]
If $ Q $ is not of the form $ \{ 1,2,\dotsc,q \} $ $ (0 \le q \le n) $, 
there exist an element $ q_0 $ in $ \widetilde{Q} $ and we fix such $ q_0 $.

\noindent
(a) For $ p \neq q_0, q_0+1 $, then 
\begin{align*}
& \sum_{ w \in \{1, w_{ q_0}^{ (p) } \} }
 w \Bigl[ \Gamma^{p,Q}(\nu) \widetilde{M}_{\nu}^{p,Q}(y) \Bigr] = 0.
\end{align*}
Here $  w_{q_0}^{ (p) }  \in {\mathcal W}_{n-1} $ permutes 
$ \nu_{q_0}^{(p)} $ and $ \nu_{q_0+1}^{(p)}  $ 
and fixes other $ \nu_i$'s.

\noindent
(b) The terms $ p = q_0 $ and $ p = q_0+1 $ cancel each other:
\begin{align*} 
& \sum_{p= q_0, q_0+1} \Gamma^{p,Q}(\nu) \widetilde{M}_{\nu}^{p,Q}(y)  = 0.
\end{align*}  
\end{itemize}
\end{prop}

\bpf 
The idea of proof is the same as Proposition \ref{claimV}.

\noindent
(i) Since $ \nu_n^{(p)} = \nu_{n-1} $ does not appear 
in the characteristic exponents of the power series $ \widetilde{M}_{\nu}^{p,Q} $, 
it is enough to show
$$
 \sum_{w \in \{ 1, w_{n-1} \} }
  w \Bigl[ \Gamma^{p,Q}(\nu) \, c_{(m_1,\dotsc,m_n)}^{p,Q}(\nu) \Bigr] = 0.
$$
By Lemma \ref{recU} (i), the recurrence relation for $ c_{(m_1,\dotsc,m_n)}^{p,Q}(\nu) $ 
is invariant under the action of $ w_{n-1} $. Actually
$ \nu_{n-1} $ does not appear
in $ \lambda_i^{p,Q} $ and 
$ \kappa^{p,Q} = \frac12 (\nu_{n-1}+\nu_n)(-\nu_{n-1}+\nu_n) $
(note that $ \widetilde{Q} = \emptyset $).
Then our task is reduced to confirm
\begin{align} \label{AA}
  \sum_{w \in \{ 1, w_{n-1} \} }
  w \Bigl[ \Gamma^{p,Q}(\nu) \, c_{(0,\dotsc,0)}^{p,Q}(\nu) \Bigr] = 0.
\end{align}
By the definition of $ \Gamma^{p,Q}(\nu) $ and Lemma \ref{recU} (ii),
\begin{align*} 
& \Gamma^{p,Q}(\nu) \, c_{(0,\dotsc,0)}^{p,Q}(\nu)
\\
& = \prod_{1 \le i<j \le n-1}
     \Gamma ( -\nu_i-\nu_j)
     \Gamma ( -\nu_i+\nu_j)
    \prod_{1 \le i \le n-1} \Gamma(-2\nu_i)
\\
& \quad \cdot
   \prod_{1 \le i \le p-1} 
     \Gamma ( -\nu_i+\nu_n)
    \prod_{p \le i \le n-1}
     \Gamma ( \nu_i-\nu_n) 
\\
& \quad \cdot
   \prod_{i \in Q} \Gamma ( -\nu_i^{(p)}-\nu_n )
   \prod_{i \in Q^c} \Gamma ( \nu_i^{(p)}+\nu_n )
 \cdot
   \frac{\Gamma (  2 \nu_n^{(p)} + 1 ) }
        {\Gamma ( \nu_n^{(p)} + \nu_n + 1 ) 
         \Gamma ( \nu_n^{(p)} - \nu_n + 1 )}.
\end{align*}
We pick up the terms containing $ \nu_n^{(p)} = \nu_{n-1} $: 
\begin{align*}
& 
\prod_{1 \le i \le n-2} 
  \Gamma ( - \nu_i^{}-\nu_{n-1})
  \Gamma ( - \nu_i^{}+\nu_{n-1})
 \cdot
  \Gamma (-2\nu_{n-1}) 
  \Gamma ( \nu_{n-1}-\nu_{n})
\\
& \quad 
\cdot 
  \Gamma ( \nu_{n-1}+\nu_n)
\cdot 
  \frac{\Gamma(2 \nu_{n-1}+1)}
       { \Gamma ( \nu_{n-1}+\nu_n+1 )
         \Gamma (\nu_{n-1}-\nu_n+1 ) }
\\
& = \frac{\pi}{ \sin(-2\nu_{n-1} \pi) }
    \cdot
    \frac{1}{(\nu_{n-1}+\nu_n)(\nu_{n-1}-\nu_n)}
    \prod_{1 \le i \le n-2} 
  \Gamma ( - \nu_i^{}-\nu_{n-1})
  \Gamma ( - \nu_i^{}+\nu_{n-1}).
\end{align*}
Then we have (\ref{AA}).
\medskip

\noindent
(ii) (a) 
In view of 
\begin{equation*}
\begin{cases}
 \nu_{q_0}^{(p)} = \nu_{q_0-1}, \ \nu_{q_0+1}^{(p)} = \nu_{q_0} & \mbox{if }  p< q_0 ,
\\
 \nu_{q_0}^{(p)} = \nu_{q_0},\ \nu_{q_0+1}^{(p)} = \nu_{q_0+1} & \mbox{if }  p> q_0+1 ,
\end{cases}
\end{equation*}
and 
$ q_0 \in Q^c $ , $ q_0+1 \in Q $, 
$ \prod_{i \in Q}(\pi y_i)^{(\cdots)} \prod_{i \in Q^c} (\pi y_i)^{(\cdots)} $ is invariant 
under the action of $ w_{q_0}^{(p)} $. As in the proof of (i), 
we can see the assertion from Lemma \ref{recU}.

%
\medskip

\noindent
(ii) (b) 
By using $ \nu_{q_0}^{(q_0)} = \nu_n $, $ \nu_{q_0+1}^{(q_0)} = \nu_{q_0} $,
$ \nu_{q_0}^{(q_0+1)} = \nu_{q_0} $ and $ \nu_{q_0+1}^{(q_0+1)} = \nu_{n} $,
our claim follows from Lemma \ref{recU}. 
Indeed we have 
\begin{align*}
& \frac{ \Gamma^{q_0,Q}(\nu) c_{(0,\dotsc,0)}^{q_0,Q}(\nu) }
  { \Gamma^{q_0+1,Q}(\nu) c_{(0,\dotsc,0)}^{q_0+1,Q}(\nu) }
\\
& = \frac{\prod_{1 \le i \le q_0-1} \Gamma(-\nu_i+\nu_n)
          \prod_{q_0 \le i \le n-1} \Gamma(\nu_i-\nu_n) }
    {\prod_{1 \le i \le q_0} \Gamma(-\nu_i+\nu_n)
          \prod_{q_0+1 \le i \le n-1} \Gamma(\nu_i-\nu_n) } 
\\
& \quad  \cdot
 \frac{ \prod_{i \in Q}\Gamma(-\nu_i^{(q_0)}-\nu_n) \prod_{i \in Q^c} \Gamma(\nu_i^{(q_0)}+\nu_n) }
 {\prod_{i \in Q}\Gamma(-\nu_i^{(q_0+1)}-\nu_n) \prod_{i \in Q^c} \Gamma(\nu_i^{(q_0+1)}+\nu_n) }
\\
& \quad \cdot \frac{ \Gamma(\nu_{q_0}^{(q_0)}-\nu_{q_0+1}^{(q_0)}+1) }
  { \Gamma(-\nu_{q_0+1}^{(q_0)} - \nu_n+1) \Gamma(\nu_{q_0}^{(q_0)}+\nu_n+1) }
  \cdot 
  \frac{\Gamma(-\nu_{q_0+1}^{(q_0+1)} - \nu_n+1) \Gamma(\nu_{q_0}^{(q_0+1)}+\nu_n+1) }
  { \Gamma(\nu_{q_0}^{(q_0+1)}-\nu_{q_0+1}^{(q_0+1)}+1) }
\\
& = \frac{ \Gamma(\nu_{q_0}-\nu_n)}{ \Gamma(-\nu_{q_0}+\nu_n) } \cdot 
    \frac{ \Gamma(-\nu_{q_0}-\nu_n) \Gamma(2\nu_n) }{ \Gamma(-2\nu_n) \Gamma(\nu_{q_0}+\nu_n) }
\\
& \quad 
    \cdot
    \frac{ \Gamma(\nu_n-\nu_{q_0}+1) }{ \Gamma(-\nu_{q_0}-\nu_n+1) \Gamma(2\nu_n+1)} \cdot 
    \frac{ \Gamma(-2\nu_n+1) \Gamma(\nu_{q_0}+\nu_n+1) }{ \Gamma(\nu_{q_0}-\nu_n+1) }
\\
& = - \frac{\sin (\nu_{q_0}+\nu_n) \pi \cdot \sin (-\nu_{q_0}+\nu_n) \pi }
    {\sin (\nu_{q_0}-\nu_n) \pi \cdot \sin (-\nu_{q_0}-\nu_n) \pi }
\\
& = -1.
\end{align*}
\epf

\bigskip 

To conclude the proof of Theorem \ref{mainthm}, we rewrite (\ref{UUU})
by using Corollary \ref{cor1} and Proposition \ref{claimU}: 
\begin{align*} 
 \widetilde{W}_{\nu}^n(y)
& =
\sum_{w \in {\mathcal W}_{n-1} }  \sum_{p=1}^{n} 
 w \biggl[ \Gamma_n(\nu^{(p)})
    \sum_{m_1,\dotsc,m_n = 0}^{\infty}  c_{n,(m_1,\dotsc, m_n)}^{}(\nu^{(p)} )
\\
& \qquad \cdot 
 \prod_{i=1}^{p-1} (\pi y_i)^{2(m_i+\nu_1+\cdots+\nu_i)}
 \prod_{i=p}^{n} (\pi y_i)^{2(m_i+\nu_1+\cdots+\nu_{i-1}+\nu_n)} \biggr]
\\
& + \sum_{w \in {\mathcal W}_{n-1} }
  \sum_{q=1}^n
  w \biggl[ \Gamma_n(\bar{\nu}^{(q)})
  \sum_{m_1,\dotsc,m_n = 0}^{\infty}  c_{n,(m_1,\dotsc, m_n)}^{}(\bar{\nu}^{(q)})
\\
& \qquad 
  \cdot \prod_{i=1}^{q-1} (\pi y_i)^{2(m_i+\nu_1+\cdots+\nu_i)}
        \prod_{i=q}^{n} (\pi y_i)^{2(m_i+\nu_1+\cdots+\nu_{i-1}-\nu_n)} \biggr]
\\
& = \sum_{ w \in {\mathcal W}_n } w \Bigl[ \Gamma_n(\nu) \widetilde{M}_{\nu}^n(y) \Bigr].
\end{align*}
\epf

\section{Mellin transforms of class one Whittaker functions}

In this section we compute the Mellin transform of the class one Whittaker function.
We first derive another recursive integral representation 
between $ W_{\nu}^n(y) $ and $ W_{\nu}^{n-1}(y) $.
This is similar to the recursive formula for $GL_n(\R)$ and $ GL_{n-2}(\R) $
Whittaker functions obtained by Stade \cite[Theorm 2.1]{St90}.
\begin{thm} \label{intrep2}
We have 
\begin{align*}
 \widetilde{W}_{\nu}^n(y)
& = 
  2^n
  \int_{(\R^{+})^{n-1}}
  \prod_{i=1}^n 
    K_{2\nu_n} \Bigl( 2\pi y_i 
   \sqrt{(1+u_{i-1}) (1+1/u_i )} \Bigr)
\\
& \qquad 
  \cdot \widetilde{W}_{\tilde{\nu}}^{n-1}
  \biggl( y_2 \sqrt{\dfrac{u_1}{u_2}}, \dotsc,
          y_{n-1} \sqrt{\dfrac{u_{n-2}}{u_{n-1}}}, y_n \sqrt{u_{n-1}} \biggr)
  \prod_{i=1}^{n-1} \dfrac{du_i}{u_i}.
\end{align*}
Here we promise $ u_0 = 1/u_n = 0 $.
\end{thm}
\bpf
Substituting $ u_i \to t_{i+1}/u_i $ $(1 \le i \le n-1)$ into the integral 
representation (\ref{intrep1}) to find 
\begin{align*}
 \widetilde{W}_{\nu}^n(y)
& = \int_{(\R^+)^n}
    \int_{(\R^+)^{n-1}}
   \prod_{i=1}^{n} 
     \exp \biggl\{ -(\pi y_i)^2 t_i - \frac{1}{t_i} \biggr\}
   \prod_{i=1}^{n} (\pi y_i t_i)^{2\nu_n}  
\\
& \qquad \cdot
   \prod_{i=1}^{n-1} 
     \exp \biggl\{ -(\pi y_i)^2 \frac{t_i}{u_i} - \frac{u_i}{t_{i+1}} \biggr\}
   \prod_{i=1}^{n-1} u_i^{-\nu_n} 
\\
& \qquad \cdot 
  \widetilde{W}_{\tilde{\nu}}^{n-1}
  \biggl( y_2 \sqrt{\frac{u_1}{u_2}}, \dotsc, 
          y_{n-1} \sqrt{\frac{u_{n-2}}{u_{n-1}}},
          y_n \sqrt{ u_{n-1} } \biggr) 
  \prod_{i=1}^{n-1} \frac{du_i}{u_i} \prod_{i=1}^{n} \frac{dt_i}{t_i} 
\\
& = \int_{(\R^+)^{n-1} }
    \prod_{i=1}^{n} \biggl( 
    \int_{0}^{\infty}
    \exp \biggl\{ -(\pi y_i)^2 t_i \Bigl(1+ \frac{1}{u_i} \Bigr)  
                 - \frac{1}{t_i} (1+u_{i-1}) \biggr\}
    \, t_i^{2\nu_n} \frac{dt_i}{t_i}  \biggr)
\\
& \qquad \cdot 
   \prod_{i=1}^{n} (\pi y_i)^{2\nu_n}  
  \prod_{i=1}^{n-1} u_i^{-\nu_n} 
  \cdot \widetilde{W}_{\tilde{\nu}}^{n-1}
  \biggl( y_2 \sqrt{\frac{u_1}{u_2}}, \dotsc, 
          y_{n-1} \sqrt{\frac{u_{n-2}}{u_{n-1}}},
          y_n \sqrt{ u_{n-1} } \biggr) 
   \prod_{i=1}^{n-1} \frac{du_i}{u_i}.
\end{align*}
By (\ref{Bessellem}), the integration $ \int dt_i $ becomes 
\begin{align*}
  2 (1+u_{i-1})^{2\nu_n} 
  \Bigl\{ (\pi y_i) \sqrt{ (1+u_{i-1})(1+1/u_i) } \Bigr\}^{-2\nu_n}
  K_{2\nu_n} \Bigl( 2\pi y_i \sqrt{(1+u_{i-1})(1+1/u_i) } \Bigr)
\end{align*}
and we finish the proof.
\epf

\bigskip

Let $ s = (s_1, \dotsc, s_n) \in \C^n $ and 
$$
 T_{\nu}^{n}(s) 
 = \int_{(\R^{+})^{n}}
   \widetilde{W}_{\nu}^{n}(y_1, \dotsc, y_n) 
    \prod_{i=1}^n (\pi y_i)^{2s_i} \frac{dy_i}{y_i} 
$$
be the multiple Mellin transform of the ($\rho$-shifted) class one Whittaker
function $ \widetilde{W}_{\nu}^n(y) $. 
In the same way as in \cite[Theorem 3.1]{St01} for $GL(n,\R)$-Whittaker functions,
we can prove the following recursive formula for $ T_{\nu}^n(s) $.
\begin{thm} Let $ n \geq 2 $ and fix real numbers $ \tau_j $ $(1 \le j \le n-1) $
such that 
$$ \tau_j 
< \min \{ {\rm Re}( \textstyle  \sum_{i=1}^j \varepsilon_i \nu_{\sigma(i)} )
             \mid \varepsilon_i \in \{ \pm 1 \},\, \sigma \in \mathfrak{S}_n \}, 
$$   
and also define
$ \tau_{-1} = + \infty $, $ \tau_0 = 0 $.
Let 
$$ \eta_j = \max \{ -\tau_{j-1} + \nu_n, \tau_{j-1} - \nu_n,
                    -\tau_j, -\tau_{j-2} \},
$$
for $ 1 \le j \le n-1 $ 
and 
$$  \Omega = 
\{ s \in \C^{n-1} \mid {\rm Re}(s_j)>\eta_j \mbox{ for } 1 \le j \le n-1 \}.
$$        
Then for $ s \in  \Omega $, we have 
\begin{align*}
 T_{\nu}^n(s) 
& = \frac{2^{-1}}{(2\pi \sqrt{-1})^{n-1}} 
  \int_{\tau_1 -\sqrt{-1} \infty}^{\tau_1+\sqrt{-1} \infty}
  \dotsi
  \int_{\tau_{n-1} -\sqrt{-1} \infty}^{\tau_{n-1}+\sqrt{-1} \infty}
  \prod_{i=1}^{n}
  \Gamma( s_i + t_{i-1} + \nu_n) \Gamma( s_i + t_{i-1} - \nu_n) 
\\
& \qquad \cdot \prod_{i=1}^{n-1}
   \frac{ \Gamma(s_i + t_i) \Gamma( s_{i+1} + t_{i-1} )  }
   { \Gamma( s_i + s_{i+1} + t_{i-1} + t_i ) }
  \cdot
   T_{\tilde{\nu}}^{n-1} (-t_1, \dotsc, -t_{n-1}) 
   \prod_{i=1}^{n-1} dt_i.
\end{align*}
Here 
$ T_{\nu}^1(s) = 2^{-1} \Gamma(s_1+\nu_1) \Gamma(s_1-\nu_1) $ 
and we promise $ t_0 = 0 $. 
\end{thm}
\bpf (sketch) 
We use induction on $n$ and 
the proof is quite analogous to \cite[\S 4]{St01}.
By Mellin inversion and Theorem \ref{intrep2},
\begin{align*}
 T_{\nu}^n(s)
&= \frac{2^{2n-1}}{(2\pi \sqrt{-1})^{n-1}}
  \int_{\tau_1 -\sqrt{-1} \infty}^{\tau_1+\sqrt{-1} \infty}
  \dotsi
  \int_{\tau_{n-1} -\sqrt{-1} \infty}^{\tau_{n-1}+\sqrt{-1} \infty} 
  T_{\tilde{\nu}}^{n-1}(-t_1,\dotsc,-t_{n-1}) 
\\
& \quad \cdot 
  \int_{(\R^{+})^{n-1}}
  \int_{(\R^+)^n}
   \prod_{i=1}^{n}
  K_{2\nu_n} \Bigl( 2 \pi y_i \sqrt{(1+u_{i-1})(1+1/u_i)} \Bigr)
  (\pi y_i)^{2s_i} 
\\
& \quad \cdot\prod_{i=1}^{n-1} 
  \biggl( \pi y_{i+1} \sqrt{\frac{u_i}{u_{i+1}}} \biggr)^{\!2t_i} 
  \prod_{i=1}^{n} \frac{dy_i}{y_i}
  \prod_{i=1}^{n-1} \frac{du_i}{u_i}
  \prod_{i=1}^{n-1} dt_i
\\
&= \frac{2^{2n-1}}{(2\pi \sqrt{-1})^{n-1}}
  \int_{\tau_1 -\sqrt{-1} \infty}^{\tau_1+\sqrt{-1} \infty}
  \dotsi
  \int_{\tau_{n-1} -\sqrt{-1} \infty}^{\tau_{n-1}+\sqrt{-1} \infty} 
  T_{\tilde{\nu}}^{n-1}(-t_1,\dotsc,-t_{n-1})   
\\
& \quad \cdot
  \int_{(\R^+)^{n-1}} 
  \prod_{i=1}^n
  \biggl(
  \int_{0}^{\infty}
  K_{2\nu_n} \Bigl( 2 \pi y_i \sqrt{(1+u_{i-1})(1+1/u_i)} \Bigr)
  (\pi y_i)^{2(s_i+t_{i-1})} \frac{dy_i}{y_i} \biggr) 
\\
& \quad \cdot 
  \prod_{i=1}^{n-1} u_i^{-t_{i-1}+t_i} 
  \prod_{i=1}^{n-1}  \frac{du_i}{u_i}
  \prod_{i=1}^{n-1} dt_i
\\
& = 
  \frac{2^{-1}}{(2\pi \sqrt{-1})^{n-1}}
  \int_{\tau_1 -\sqrt{-1} \infty}^{\tau_1+\sqrt{-1} \infty}
  \dotsi
  \int_{\tau_{n-1} -\sqrt{-1} \infty}^{\tau_{n-1}+\sqrt{-1} \infty} 
  T_{\tilde{\nu}}^{n-1}(-t_1,\dotsc,-t_{n-1})
\\
& \quad \cdot 
  \prod_{i=1}^{n}
  \Gamma( s_i+t_{i-1}+\nu_n) \Gamma(s_i+t_{i-1}-\nu_n) 
\\
& \quad \cdot 
  \int_{(\R^+)^{n-1}}
  \biggl\{ 
  \prod_{i=1}^{n-1} u_i^{-t_{i-1}+t_i} 
  \prod_{i=1}^{n}
  \{(1+u_{i-1})(1+1/u_i)\}^{-s_i-t_{i-1}} \biggr\}
  \prod_{i=1}^{n-1} \frac{du_i}{u_i}
\prod_{i=1}^{n-1} dt_i
\\
& = 
\frac{2^{-1}}{(2\pi \sqrt{-1})^{n-1}}
  \int_{\tau_1 -\sqrt{-1} \infty}^{\tau_1+\sqrt{-1} \infty}
  \dotsi
  \int_{\tau_{n-1} -\sqrt{-1} \infty}^{\tau_{n-1}+\sqrt{-1} \infty} 
  T_{\tilde{\nu}}^{n-1}(-t_1,\dotsc,-t_{n-1})
\\
& \quad \cdot 
  \prod_{i=1}^{n}
  \Gamma( s_i+t_{i-1}+\nu_n) \Gamma(s_i+t_{i-1}-\nu_n) 
\\
& \quad \cdot 
 \prod_{i=1}^{n-1}
 \biggl( \int_0^{\infty} (1+u_i)^{-s_i-s_{i+1}-t_{i-1}-t_i}
   u_i^{s_i+t_i} \frac{du_i}{u_i} \biggr)
  \prod_{i=1}^{n-1} dt_i.
\end{align*}
By using 
$$
 \int_0^{\infty}
  (1+u)^{-(x+y)} u^{y} \frac{du}{u}
= \frac{ \Gamma(x) \Gamma(y) }{ \Gamma(x+y)}
$$
for $ {\rm Re}(x)>0 $, $ {\rm Re}(y)>0 $, we get the assertion.
$ \Box $



\begin{thebibliography}{99}  
\bibitem{Bu}
Bump, D., 
Automorphic forms on $GL(3,{\mathbf R})$, 
Lect. Note in Math. {\bf 1083} (1984), Springer-Verlag. 


 
\bibitem{Ha0}
M. Hashizume,
Whittaker models for real reductive groups,
Japan. J. Math. {\bf 5} (1979), 349--401.

\bibitem{Ha}
M. Hashizume, 
Whittaker functions on semisimple Lie groups,
Hiroshima Math. J. {\bf 12} (1982), 259--293. 

\bibitem{HC}
Harish-Chandra,
Spherical functions on a semisimple Lie group I,II,
Amer. J. Math. {\bf 80} (1958), 241--310, 553--613.

\bibitem{HIO}
M. Hirano, T. Ishii and T. Oda,   
Whittaker functions for $P_J$-principal series representations of $Sp(3,\R)$,
Adv. Math. {\bf 215} (2007), 734--765.
   
\bibitem{I1}
T. Ishii,
A remark on Whittaker functions on $SL(n,\R)$,
Ann. Inst. Fourier {\bf 55} (2005), 483--492.

\bibitem{I2}
T. Ishii,
Principal series Whittaker functions on $Sp(2,\R)$,
J. Funct. Anal. {\bf 225} (2005), 1--32.

\bibitem{I3}
T. Ishii,
Class one Whittaker functions on real semisimple Lie groups,
Automorphic representations, $L$-functions and periods,
S\=urikaisekikenky\=usho K\=oky\=uroku No. {\bf 1523} (2006), 70--78.


\bibitem{ISt}
T. Ishii and E. Stade,
New formulas for Whittaker functions on $ GL(n,\R)$, 
J. Funct. Anal. {\bf 244} (2007), 289--314.

\bibitem{Ja}
H. Jacquet,
Fonctions de Whittaker associ\'ees aux groupes de Chevalley,
Bull. Soc. Math. France {\bf 95} (1967), 243--309.


\bibitem{Ko}
B. Kostant,
On Whittaker vectors and representation theory,
Invent. Math. {\bf 48} (1978), 101--184.

\bibitem{Sh}
J. Shalika,
The multiplicity one theorem for $GL(n)$,
Ann. of Math. {\bf 100} (1974), 171--193.


\bibitem{St90}
E. Stade, 
On explicit integral formulas for $GL(n,\R)$-Whittaker functions,
Duke Math. J. {\bf 60} (1990), 313--362.

\bibitem{St93}
E. Stade,
$GL(4,{\mathbf R})$ Whittaker functions and ${}_4F_3(1)$ hypergeometric series, 
Trans. Amer. Math. Soc.  {\bf 336} (1993), 253--264. 

\bibitem{St01}
E. Stade,
Mellin transforms of $GL(n,\R)$ Whittaker functions,
Amer. J. Math. {\bf 123} (2001), 121--161.

\bibitem{St02}
E. Stade,
Archimedean $L$-factors on $GL(n) \times GL(n)$ and generalized Barnes integrals,
Israel J. Math. {\bf 127} (2002), 201-220.



\bibitem{Wa}
N. Wallach,
Asymptotic expansions of generalized matrix entries of representations of 
real reductive groups,
Lect. Notes in Math. {\bf 1024} (1984), 287--369.
\end{thebibliography}
\end{document}